\documentclass{amsart}
\usepackage{amssymb, amsbsy, amsthm, amsmath,amstext, amsopn}
\usepackage[all]{xy}
\usepackage{amsfonts}
\usepackage{amscd}
\usepackage{mathrsfs}
\usepackage{verbatim}

\newtheorem{thm}{Theorem} [section]
\newtheorem{lemma}[thm]{Lemma}

\newtheorem{corollary}[thm]{Corollary}
\newtheorem{prop}[thm]{Proposition}
\newtheorem{notation}[thm]{Notation}

\theoremstyle{definition}

\newtheorem{defn}[thm]{Definition}

\newtheorem{example}[thm]{Example}

\theoremstyle{remark}

\newtheorem{remark}[thm]{Remark}

\begin{document}

\numberwithin{equation}{section}

\newcommand{\hs}{\mbox{\hspace{.4em}}}
\newcommand{\ds}{\displaystyle}
\newcommand{\bd}{\begin{displaymath}}
\newcommand{\ed}{\end{displaymath}}
\newcommand{\bcd}{\begin{CD}}
\newcommand{\ecd}{\end{CD}}

\newcommand{\on}{\operatorname}
\newcommand{\proj}{\operatorname{Proj}}
\newcommand{\bproj}{\underline{\operatorname{Proj}}}
\newcommand{\spec}{\operatorname{Spec}}
\newcommand{\bspec}{\underline{\operatorname{Spec}}}
\newcommand{\pline}{{\mathbf P} ^1}
\newcommand{\aline}{{\mathbf A} ^1}
\newcommand{\pplane}{{\mathbf P}^2}
\newcommand{\aplane}{{\mathbf A}^2}
\newcommand{\coker}{{\operatorname{coker}}}
\newcommand{\ldb}{[[}
\newcommand{\rdb}{]]}

\newcommand{\Sym}{\operatorname{Sym}^{\bullet}}
\newcommand{\Symp}{\operatorname{Sym}}
\newcommand{\Pic}{\operatorname{Pic}}
\newcommand{\AAut}{\operatorname{Aut}}
\newcommand{\PAut}{\operatorname{PAut}}

\newcommand{\bc}{{\bf BC}\, }
\newcommand{\ec}{{\bf EC}\,}

\newcommand{\ul}{\underline}
\newcommand{\too}{\twoheadrightarrow}
\newcommand{\C}{{\mathbb C}}
\newcommand{\Z}{{\mathbf Z}}
\newcommand{\R}{{\mathbf R}}
\newcommand{\Cx}{{\mathbf C}^{\times}}
\newcommand{\Cbar}{\overline{\C}}
\newcommand{\Cxbar}{\overline{\Cx}}
\newcommand{\cA}{{\mathcal A}}
\newcommand{\cS}{{\mathcal S}}
\newcommand{\cV}{{\mathcal V}}
\newcommand{\cM}{{\mathcal M}}
\newcommand{\bA}{{\mathbf A}}
\newcommand{\cB}{{\mathcal B}}
\newcommand{\cC}{{\mathcal C}}
\newcommand{\cD}{{\mathcal D}}
\newcommand{\D}{{\mathcal D}}
\newcommand{\cs}{{\mathbf C} ^*}
\newcommand{\boldc}{{\mathbf C}}
\newcommand{\cE}{{\mathcal E}}
\newcommand{\cF}{{\mathcal F}}
\newcommand{\bF}{{\mathbf F}}
\newcommand{\cG}{{\mathcal G}}
\newcommand{\G}{{\mathbb G}}
\newcommand{\cH}{{\mathcal H}}
\newcommand{\cJ}{{\mathcal J}}
\newcommand{\cK}{{\mathcal K}}
\newcommand{\cL}{{\mathcal L}}
\newcommand{\baL}{{\overline{\mathcal L}}}
\newcommand{\bL}{{\mathbf L}}
\newcommand{\M}{{\mathcal M}}
\newcommand{\Mf}{{\mathfrak M}}
\newcommand{\bM}{{\mathbf M}}
\newcommand{\bm}{{\mathbf m}}
\newcommand{\cN}{{\mathcal N}}
\newcommand{\theo}{\mathcal{O}}
\newcommand{\cP}{{\mathcal P}}
\newcommand{\cT}{{\mathcal T}}
\newcommand{\cR}{{\mathcal R}}
\newcommand{\Pp}{{\mathbb P}}
\newcommand{\boldp}{{\mathbf P}}
\newcommand{\boldq}{{\mathbf Q}}
\newcommand{\bbL}{{\mathbf L}}
\newcommand{\cQ}{{\mathcal Q}}
\newcommand{\cO}{{\mathcal O}}
\newcommand{\Oo}{{\mathcal O}}
\newcommand{\OX}{{\Oo_X}}
\newcommand{\OY}{{\Oo_Y}}
\newcommand{\otY}{{\underset{\OY}{\ot}}}
\newcommand{\otX}{{\underset{\OX}{\ot}}}
\newcommand{\cU}{{\mathcal U}}\newcommand{\cX}{{\mathcal X}}
\newcommand{\cW}{{\mathcal W}}
\newcommand{\boldz}{{\mathbf Z}}
\newcommand{\qgr}{\operatorname{qgr}}
\newcommand{\gr}{\operatorname{gr}}
\newcommand{\rk}{\operatorname{rk}}
\newcommand{\coh}{\operatorname{coh}}
\newcommand{\End}{\operatorname{End}}
\newcommand{\uEnd}{\underline{\operatorname{End}}}
\newcommand{\Hom}{\operatorname{Hom}}
\newcommand{\uHom}{\underline{\operatorname{Hom}}}
\newcommand{\uHomY}{\uHom_{\OY}}
\newcommand{\uHomX}{\uHom_{\OX}}
\newcommand{\Ext}{\operatorname{Ext}}
\newcommand{\bExt}{\operatorname{\bf{Ext}}}
\newcommand{\Tor}{\operatorname{Tor}}

\newcommand{\inv}{^{-1}}
\newcommand{\airtilde}{\widetilde{\hspace{.5em}}}
\newcommand{\airhat}{\widehat{\hspace{.5em}}}
\newcommand{\nt}{^{\circ}}
\newcommand{\del}{\partial}
\newcommand{\Qgr}{\on{Qgr}}

\newcommand{\supp}{\operatorname{supp}}
\newcommand{\GK}{\operatorname{GK-dim}}
\newcommand{\hd}{\operatorname{hd}}
\newcommand{\id}{\operatorname{id}}
\newcommand{\res}{\operatorname{res}}
\newcommand{\lrar}{\leadsto}
\newcommand{\im}{\operatorname{Im}}
\newcommand{\HH}{\operatorname{H}}
\newcommand{\TF}{\operatorname{TF}}
\newcommand{\Bun}{\operatorname{Bun}}
\newcommand{\BunD}{\operatorname{Bun}_{\D^{\bullet}}}
\newcommand{\PicD}{\operatorname{Pic}_{\D}}
\newcommand{\Hilb}{\operatorname{Hilb}}
\newcommand{\Fact}{\operatorname{Fact}}
\newcommand{\CM}{\mathsf{CM}}
\newcommand{\PB}{\mathsf{PB}}
\newcommand{\Kos}{\mathsf{Kos}}
\newcommand{\mOp}{\mathcal{MO}p}
\newcommand{\MD}{\mathfrak{M}^{\D}}
\newcommand{\F}{\mathcal{F}}
\newcommand{\Ff}{\mathbb{F}}
\newcommand{\nthord}{^{(n)}}
\newcommand{\Aut}{\underline{\operatorname{Aut}}}
\newcommand{\Gr}{{\mathfrak{Gr}}}
\newcommand{\GR}{\on{Gr}}
\newcommand{\GRo}{{\mathfrak{GR}^{\circ}}}
\newcommand{\Fr}{\operatorname{Fr}}
\newcommand{\GL}{\operatorname{GL}}
\newcommand{\gl}{\mathfrak{gl}}
\newcommand{\SL}{\operatorname{SL}}
\newcommand{\ff}{\footnote}
\newcommand{\ot}{\otimes}
\def\Ext{\operatorname {Ext}}
\def\Hom{\operatorname {Hom}}
\def\Ind{\operatorname {Ind}}
\def\bbZ{{\mathbb Z}}

\newcommand{\AutO}{\on{Aut}\Oo}
\newcommand{\Der}{{\on{Der}\,}}
\newcommand{\DerO}{{\on{Der}\Oo}}
\newcommand{\AutK}{\on{Aut}\K}

\newcommand{\nc}{\newcommand}
\nc{\ol}{\overline}
\nc{\cont}{\on{cont}}
\nc{\rmod}{\on{mod}}
\nc{\Mtil}{\widetilde{M}}
\nc{\wb}{\overline}
\nc{\wt}{\widetilde}
\nc{\wh}{\widehat}
\nc{\sm}{\setminus}
\nc{\mc}{\mathcal}
\nc{\mbb}{\mathbb}
\nc{\Mbar}{\wb{M}}
\nc{\Nbar}{\wb{N}}
\nc{\Mhat}{\wh{M}}
\nc{\pihat}{\wh{\pi}}
\nc{\JYX}{\cJ_{Y\leftarrow X}}
\nc{\phitil}{\wt{\phi}}
\nc{\Qbar}{\wb{Q}}
\nc{\DYX}{\D_{Y\leftarrow X}}
\nc{\DXY}{\D_{X\to Y}}
\nc{\dR}{\stackrel{\bbL}{\underset{\D_X}{\ot}}}
\nc{\Winfi}{\cW_{1+\infty}}
\nc{\K}{{\mc K}} \nc{\Kx}{{\mc K}^{\times}}
\nc{\Ox}{{\mc O}^{\times}}
\nc{\unit}{{\bf \on{unit}}}
\nc{\boxt}{\boxtimes}
\nc{\xarr}{\stackrel{\rightarrow}{x}}

\nc{\Gamx}{\Gamma_-^{\times}}
\nc{\Gap}{\Gamma_+}
\nc{\Emx}{\cE_-^{\times}}

\nc{\Coh}{{\mathcal Coh}} \nc{\Qcoh}{{\mathcal Qcoh}}
\nc{\AdR}{A_{dR}^\cdot} \nc{\Adol}{A_{Dol}^\cdot} \nc{\ST}{\Sym\cT}
\nc{\Conn}{{\mc Conn}}

\nc{\Enat}{E^{\natural}}
\nc{\Enatbar}{\ol{E}^{\natural}}
\nc{\Of}{{\mathbb O}}
\nc{\obar}{\wb{o}}
\newcommand{\SShv}{\on{SShv}}

\nc{\Ga}{\G_a}

\nc{\Gm}{\G_m} \nc{\Gabar}{\wb{\G}_a} \nc{\Gmbar}{\wb{\G}_m}
\nc{\PD}{{\mathbb P}_{\D}} \nc{\Pbul}{P_{\bullet}}
\nc{\Tors}{\on{Tors}} \nc{\PS}{{\mathsf{PS}}}

\title{Perverse Bundles and Calogero-Moser Spaces}
\author{David Ben-Zvi}
\address{Department of Mathematics\\University of Texas\\Austin, TX
78712-0257}
\email{benzvi@math.utexas.edu}
\author{Thomas Nevins}
\address{Department of Mathematics\\University of Illinois at
Urbana-Champaign\\Urbana, IL 61801}
\email{nevins@uiuc.edu}
\begin{abstract}
We present a simple description of moduli spaces of torsion-free
$\D$-modules ($\D$-{\em bundles}) on general smooth complex curves,
generalizing the identification of the space of ideals in the Weyl
algebra with Calogero-Moser quiver varieties. Namely, we show that
the moduli of $\D$-bundles form twisted cotangent bundles to moduli
of torsion sheaves on $X$, answering a question of Ginzburg. The
corresponding (untwisted) cotangent bundles are identified with
moduli of {\em perverse vector bundles} on $T^*X$, which contain as
open subsets the moduli of framed torsion-free sheaves (the Hilbert
schemes $T^*X^{[n]}$ in the rank one case). The proof is based on
the description of the derived category of $\D$-modules on $X$ by a
noncommutative version of the Beilinson transform on $\pline$.
\end{abstract}

\maketitle
\section{Introduction}
The starting point for this paper is the following.  Let
\bd
M = \gl_n\times \gl_n \times {\mathbb C}^n\times ({\mathbb C}^n)^* =
T^*(\gl_n\times {\mathbb C}^n).\ed
  Define a map $\mu: M\rightarrow \gl_n$
by $\mu(X,Y,i,j) = [X,Y] + ij$. It is well known (see e.g.
\cite{Nakajima}) that the geometric invariant theory quotient (at a
nontrivial character $\chi$, the determinant character of $GL_n$)
$$(\C^2)^{[n]}=\mu^{-1}(0)/\!\!/_{\chi}GL_n$$ is the Hilbert scheme of $n$
points on $\C^2$. Indeed,
stability forces $j$ to vanish and $i\in\C^n$ to generate $\C^n$
under multiplication by $X$ and $Y$, so $\mu^{-1}(0)$ parametrizes
commuting $n\times n$ matrices $X,Y$
 with a cyclic vector for $\C[X,Y]$; hence $\mu^{-1}(0)/\!\!/GL_n$ parametrizes length $n$ quotients of $\C[x,y]$.

Now consider the spaces \bd \CM_n = \mu^{-1}(I)/GL_n, \ed known as
the {\em $n$-particle rational Calogero-Moser (CM) spaces}. The CM
spaces arise in the study of integrable systems and soliton
equations \cite{Wilson CM}; they
also play a central role in the representation theory of Cherednik
algebras \cite{EG}. We write $A_1=\D_{\aline}= {\mathbb
C}[x,\partial/\partial x]$ for the {\em first Weyl algebra}. There
is an appealing generalization by Berest and Wilson \cite{Wilson CM,
BW automorphisms, BW ideals} (following earlier work of
Cannings-Holland \cite{CH ideals} and Le Bruyn \cite{Le Bruyn}) of
the above description of the Hilbert scheme to the classification of
ideals in $A_1$:
\begin{thm}[\cite{BW automorphisms,BW ideals}]\label{the BW thm}
The space $\CM_n$ parametrizes (isomorphism classes of) right
ideals in $A_1$ that ``have second Chern class $n$.''
\end{thm}
It is tempting to think that the meaning of the theorem is
transparent, since the CM relation $\mu(X,Y,i,j)=I$, that is, $[X,Y]
+ ij = I$, is very close to the defining relation
$[x,\partial/\partial x] = 1$ for $A_1$.  However, the equation
$[X,Y]=I$ for $n\times n$ matrices has no solutions, so the naive
generalization of the description of the Hilbert scheme fails.

The
relations between Hilbert schemes, Calogero-Moser spaces and the
Weyl algebra have been further explored in many works.  These
include \cite{KKO,BGK1,BGK2, solitons}, and, most notably for the
present paper, \cite{BC}, where $\CM_n$ was described via
$A_\infty$-modules over $A_1$. In \cite{KKO}, the moduli spaces of
({\em framed}) modules of any rank over $A_1$ were described and
related to the classical (ADHM quiver) description of moduli of
torsion-free sheaves on $\C^2$ (framed at infinity), the
Berest-Wilson theorem being the case of rank one; in \cite{BGK1, BGK2}, the
corresponding moduli spaces for a more general class of
``noncommutative planes'' were identified with quiver varieties.

In light of the importance of the CM spaces, one naturally wants to
generalize the description of ideals in $A_1=\D_{\aline}$ to that of
torsion-free modules over differential operators on a higher genus
curve (the special case of elliptic curves was treated in
\cite{solitons} by methods that do not generalize to higher genus).
As we indicated already, the space $M$ with which we started is the
cotangent bundle of $Q_n = \gl_n\times {\mathbb C}^n$. Furthermore,
$Q_n$ is identified with a type of Quot scheme: it is the moduli
space of quotients $\theo_{{\mathbf A}^1}^n\rightarrow Q$, where $Q$
is a torsion sheaf on the affine line of length $n$ equipped with a
section $i\in H^0(Q)$.  Finally, $\mu$ is identified with a moment
map for the $GL_n$-action on $M=T^*Q_n$. With this in mind, Ginzburg
proposed:
\begin{defn}[Ginzburg]
Let $X$ be a smooth curve and $Q_n(X)$ the variety parametrizing pairs
$(\theo^n_X\rightarrow Q, i\in H^0(Q))$ where $Q$ is a length $n$ torsion sheaf
on $X$.     The {\em $n$th CM space of $X$}, denoted $\CM_n(X)$,
is the Hamiltonian reduction of $T^*Q_n(X)$ at $I\in \gl_n=\gl_n^*$.
\end{defn}

The following theorem answers a question of Ginzburg asking for the
higher-genus generalization of Theorem \ref{the BW thm}:
\begin{thm}\label{state main thm}
The space $\CM_n(X)$ is the ``Hilbert scheme of the quantized
cotangent bundle of $X$," i.e. the moduli space for trivially-framed
$\D$-line bundles (Definition \ref{D-bundle def} and Proposition
\ref{D-bundle equiv}) on $X$ with second Chern class $n$.
\end{thm}
\noindent Our method, which we outline below, gives as a byproduct a
quite transparent description of the procedures that identify rank
$1$ torsion-free $A_1$-modules with quadruples $(X,Y,i,j)$
satisfying the CM relation as an application of Koszul duality.

In fact, we prove a significantly more general result. Let $V$
denote a fixed vector bundle on $X$.  Let $\PS(X,V)$ denote the {\em
perverse symmetric power} of the curve $X$; this is the moduli stack
of framed torsion sheaves on $X$. i.e. pairs $(Q,i)$, where $Q$ is a
torsion coherent sheaf on $X$ and $i:V\to Q$ is a homomorphism. The
substack where the length of $Q$ is $n$ is denoted $\PS_n(X,V)$. The
terminology is motivated by the embedding of the symmetric powers
${\mathsf{Sym}^n X}\hookrightarrow \PS(X,\Oo)$ as the locus where
$i$ is surjective (i.e. where the complex $\Oo\to Q$ has no first
cohomology).

A
 {\em perverse vector bundle} on the compactification
$S={\mathbb P}(T_X\oplus \theo_X)$ of $T^*X$ is
 a coherent complex $\cF$ on $S$ with $\HH^0(\cF)$
torsion-free, $\HH^1(\cF)$ zero-dimensional and all other cohomology
sheaves vanishing (Definition \ref{D-bundle def}); it is {\em
$V$-framed} if its restriction to the divisor ${\mathbb P}(T_X\oplus
\theo_X)\setminus T^*X \cong X$ is equipped with a quasi-isomorphism
to $V$. A {\em $\D$-bundle on $X$} is a torsion-free $\D$-module on $X$, which
one may consider as a torsion-free sheaf on the quantized cotangent
bundle to $X$. A $V$-framing on a $\D$-bundle is the noncommutative
analog of fixing the behavior of a torsion-free sheaf on $T^*X$ near
the curve at infinity in $S$. Technically it amounts to a filtration
of a prescribed form---a precise definition is given in
Definition \ref{D-bundle def} below (which, by
Proposition \ref{D-bundle equiv}, agrees with our previous definition in
\cite{solitons}).

\begin{thm}\label{state main theorem}
\mbox{}
\begin{enumerate}
\item The cotangent bundle of the perverse symmetric power $\PS(X,V)$ is the
moduli stack $\PB(X,V)$ of $V$--framed perverse vector bundles on
$T^*X$.

\item The {\em Calogero-Moser space} $\CM(X,V)$, defined as the
twisted cotangent bundle (see \cite{BB})
 of $\PS(X,V)$ associated to the dual
determinant line bundle, is isomorphic to the moduli stack of
$V$-framed $\D$-bundles on $X$.
\end{enumerate}
Moreover, the components of $\CM$ and $\PB$ over $\PS_n(X,V)$
parametrize $\D$-bundles (respectively perverse bundles) with {\em
second Chern class} $c_2=n$.
\end{thm}

The reader may wish to view Theorem \ref{state main theorem} in
analogy with a well-known description of the moduli of flat vector
bundles on $X$. Let $\Bun_n X$ denote the moduli stack of rank $n$
vector bundles on $X$. The cotangent bundle $T^*\Bun_n X$ is
identified with the moduli stack of Higgs bundles on $X$, which are
coherent sheaves on $T^*X$ which are finite of degree $n$ over $X$.
On the other hand, the moduli space of rank $n$ bundles with flat
connection is the twisted cotangent bundle of $\Bun_n X$
corresponding to the dual of the determinant line bundle, i.e. the
space of connections on this bundle (see for example
\cite{Faltings,szego}). This twisted cotangent bundle is not
``supported'' everywhere: for a bundle to admit a connection, all its
indecomposable summands must have degree zero, so the image in
$\Bun_n(X)$ of this twisted cotangent bundle
 is a proper subset of $\Bun_n(X)$.   However, it carries a
canonical action of the cotangent bundle and is a torsor for
$T^*\Bun_nX$ over its image in $\Bun_n(X)$---which we refer to as
its {\em support}---making it a {\em pseudo-torsor} over its
support.

Let us spell out more explicitly the meaning of part (2) of Theorem
\ref{state main theorem}. Let $\wt{\PS}(X,V)$ denote the (Quot-type)
scheme parametrizing $(Q,i)$ as above together with an
identification $\Gamma(Q)={\mathbb C}^n$ (this scheme is actually a smooth
variety). The Calogero-Moser space
$\CM(X,V)$ is given by hamiltonian reduction by $GL_n$ of the
cotangent bundle of $\wt{\PS}(X,V)$, with moment map given by the
dual to the determinant character $\on{det}\in\gl_n^*$. This is a
pseudo-torsor for the cotangent bundle of $\PS(X,V)$, i.e. a torsor
over the locus of its support. One can show that the support of
$\CM(X,V)\to {\PS}(X,V)$ is the substack of indecomposable framed
torsion sheaves.

The moduli stack $\PB_n(X,\Oo^k)$ contains the moduli of framed rank
$k$ torsion-free sheaves on $T^*X$, as the open subset of perverse
bundles with vanishing first cohomology. In particular, in the rank
one case we find that the Hilbert scheme of $n$ points on $T^*X$ is
an open subset of the moduli of the {\em perverse Hilbert scheme}:
$$(T^*X)^{[n]}\subset \PB_n(X)=T^*\PS_n(X)$$
(where we drop the trivial framing from the notation). In the case
of $\aline$, the perverse Hilbert scheme is the Hamiltonian
reduction $\mu\inv(0)/GL_n$, i.e. we have an open embedding
$$(\C^2)^{[n]}\subset \PB_n(\aline)=T^*[(\gl_n\times\C^n)/GL_n]$$
obtained by dropping the stability condition for the GIT quotient.

On the other hand, Theorem \ref{state main theorem} asserts in
particular that the moduli of (trivially framed $c_2=n$)
$\D$-bundles is the pseudo-torsor $\CM_n(X)$ over the perverse
Hilbert scheme $\PB_n(X)=T^*\PS_n(X)$. Moreover, if $\on{Pic}(X)$ is
trivial, rank one $\theo_X$-framed $\D$-bundles up to isomorphism
correspond to isomorphism classes of ideals in $\D$. For $X=\aline$,
we thus recover the intimate relationship between Calogero-Moser
spaces, ideals in the Weyl algebra and the Hilbert scheme of points
in the plane that appears in the numerous works cited above.

\subsection{Techniques and Outline}
We will deduce Theorem \ref{state main theorem} from a general
description of arbitrary framed $\D$-modules on (not necessarily
projective) smooth curves $X$ (in particular, of general holonomic
$\D$-modules and $\D$-bundles), generalizing a classical description
of flat vector bundles (which are precisely $\D$-modules framed by
$0$). As we review in Section \ref{NC background}, framed
$\D$-modules on $X$ and framed sheaves on $T^*X$ are both examples
of sheaves on noncommutative $\pline$-bundles (i.e. noncommutative
ruled surfaces) over $X$, in the sense of \cite{VdB2}. In order to
describe sheaves on a $\pline$-bundle, it is natural to imitate the
Beilinson transform description of complexes on projective spaces,
or equivalently to apply Koszul duality. This technique was first
applied to give a description of $\D$-bundles on projective curves
by Katzarkov, Orlov and Pantev \cite{KOP}. We build on this idea to
prove a general derived equivalence (Theorem \ref{state general
framed theorem}) from which we deduce descriptions of various
classes of objects (Theorems \ref{CM thm} and  \ref{ruled surface
thm}).

As an important technical note, we work throughout not with derived
categories of modules as triangulated categories, but rather with
differential graded (dg) categories which are ``enhancements" of
underlying (triangulated) derived categories. See \cite{Keller} for
an excellent overview and Section \ref{dg derived cats} for more
discussion of the framework we need.

In Section \ref{NC background} we develop the general algebraic
setting we need for noncommutative $\pline$-bundles, by localizing
the category of graded modules for the Rees algebra $\cR(\D)$ with
respect to bounded modules (the usual construction of
$\Qgr\,\cR(\D)$).  We also introduce the dg derived category,
denoted $D_{dg}(\PD)$, of this localized module category. In Section
\ref{Beilinson}, we adapt the Beilinson transform to our setting and
develop an analog of \v{C}ech cohomology for computing its output.
We then prove, in Section \ref{identifying moduli}, that the
Beilinson transform gives a concise description of $D_{dg}(\PD)$: it
is equivalently described as a dg category of {\em Koszul data}
(Section \ref{defn of Kos}).\footnote{Alternatively, one may use the
Koszul duality between modules for the Rees algebra and for the de
Rham algebra \cite{Kap} to describe $D_{dg}(\PD)$.}    An object
$M\in D_{dg}(\PD)$ has a well-defined ``restriction to the curve at
infinity" $i_\infty^* M$, and we may specialize the general Koszul
description to framed complexes:

\begin{thm}\label{state general framed theorem}(Theorem \ref{general
framed thm}) There is a natural quasi-equivalence of dg categories
between $D_{dg}(\PD)$ and the dg category $\Kos$ of Koszul data (see
Section \ref{defn of Kos}) whose objects are triples $\cC = (C_{-1},
C_0, a: \D^1\otimes C_{-1}\rightarrow C_0)$ consisting of objects
$C_{-1}$ and $C_0$ of $D_{dg,\on{qcoh}}(X)$ together with a morphism
$a:\D^1\otimes C_{-1}\rightarrow C_0$ of complexes.

Under this quasi-equivalence, a choice of $V$-framing of an object
$M$ of $D_{dg}(\PD)$ corresponds to a choice of quasi-isomorphism
$\on{Cone}(C_{-1} = \theo\otimes C_{-1} \xrightarrow{a} C_0)\simeq
V$ for the corresponding object $(C_{-1},C_0,a)$ of $\Kos$.
\end{thm}
An identical result holds in the commutative case, i.e. for
complexes on $\overline{T^*X}$, with $\D^1$ replaced by $\Oo\oplus
\cT$: this is just the families version of the usual description of
the derived category of $\pline$, and the same methods can be
extended to the noncommutative setting. Furthermore, the
 difference of two $\D^1$-action maps $a_1,a_2$ as above gives
their commutative analog, thereby explicitly exhibiting a
pseudotorsor structure for framed $\D$-complexes over framed
complexes on $T^*X$. This is one of the main points of our paper:
that the standard Beilinson description of sheaves on the projective
line, with straightforward modifications, provides an explicit description of
the moduli of $\D$-modules on curves generalizing and clarifying the
much-studied case of the Weyl algebra.

Next, in Section \ref{identifying moduli}, we identify which objects
in $\Kos$ correspond to honest framed $\D$-modules, rather than
complexes. Derived categories of surfaces carry a natural perverse
coherent $t$-structure, obtained from the standard $t$-structure by
tilting the torsion sheaves of dimension zero into cohomological
degree one.\footnote{A $t$-structure on a dg category
is just a $t$-structure on its derived category (Section \ref{dg derived cats}).}
 A similar definition makes sense for the noncommutative
ruled surface defined by $\D$-modules; however, since there are no
$\D$-modules with zero-dimensional support, this recovers the
standard $t$-structure. As might be expected from the interpretation
as a Koszul duality or de Rham functor, it is the perverse
$t$-structure which is compatible with the Beilinson transform
description above. In particular we find that the commutative analog
of the data describing $\D$-modules parametrizes (framed) perverse
coherent sheaves on $T^*X$.

The equivalence of Theorem \ref{state general framed theorem} simplifies
considerably in the case of $\D$-bundles, i.e., the case of pure two-dimensional
support.  Suppose the framing $V$ is a vector bundle.  Then, under
the quasi-equivalence of Theorem \ref{state general framed theorem},
a $V$-framed perverse $\D$-bundle corresponds to a triple
$(C_{-1}, C_0, a)$ in which  $C_{-1}=Q[-1]$ where
$Q$ is a torsion sheaf on $X$, and $C_0=\on{Cone}(V\xrightarrow{i}
Q)$ for a map $i$. We thus obtain:
\begin{thm}[Theorem \ref{perverse thm}]\label{CM thm}
For a vector bundle $V$ on $X$, the moduli stack $\CM(X,V)$ is
isomorphic to the stack of triples $(Q,i, a)$ where $Q$ is a torsion sheaf,
$i:V\rightarrow Q$ is a homomorphism, and
\bd
a: Q[-1]\rightarrow \on{Cone}(V\xrightarrow{i} Q)
\ed
is a map in $D_{dg, qcoh}(X)$.
\end{thm}
\noindent
The corresponding commutative data give
perverse vector bundles on $T^*X$.  More generally,
let
\bd
0\rightarrow \theo_X\rightarrow {\mathcal E} \rightarrow \cE/\theo_X\rightarrow 0
\ed
be an exact sequence of vector bundles on $X$, with $\on{rk}(\cE) =2$.
Let $S={\mathbf P}(\cE)$ with the section $\sigma$ determined by $\theo_X\subset \cE$.
\begin{thm}\label{ruled surface thm}
\mbox{}
\begin{enumerate}
\item The moduli stack
of $V$-framed perverse vector bundles on $S$ (framed along the section $\sigma$) is isomorphic to the moduli
stack of triples $(Q,i,a)$ where $(Q,i)$ is a framed torsion sheaf and
\bd
a:\cE\otimes Q[-1]\rightarrow \on{Cone}(V\xrightarrow{i}Q)
\ed
is a map whose restriction to $\theo\otimes Q[-1]$ is the natural "inclusion."
\item In particular, the conclusion of Theorem \ref{perverse thm} holds for
$\PB(X,V)$, with $\D^1$ replaced by $\Oo\oplus{\mc T}_X$.
\end{enumerate}
\end{thm}

\begin{remark}[Holonomic modules] It is easy to read off from Theorem
\ref{state general framed theorem} the analog of the above
classification in the case of holonomic $\D$-modules, i.e. the case
of pure one-dimensional support.  More precisely, for a torsion
sheaf $V$ on $X$, the moduli stack of $V$-framed (hence holonomic)
$\D$-modules is isomorphic to the stack of triples $(C_{-1},C_0, a)$ as in
Theorem \ref{state general framed theorem}, where $C_{-1}$ is a
vector bundle (in degree zero) and $C_0$ is an extension of $V$ by $C_{-1}$.
 The same data, with $\D^1$
replaced by $\Oo\oplus{\mc T}_X$, describes the moduli stack of
$V$-framed coherent sheaves on $T^*X$ with pure one-dimensional
support---these are precisely the {\em framed spectral sheaves} on
$T^*X$, as described in \cite{KOP,flows}.
\end{remark}

 The
resulting descriptions are evidently functorial and local in $X$ (in
other words, we describe framed $\D$-modules and perverse coherent
sheaves as stacks over $X$). When the curve is projective, Serre
duality then identifies the perverse coherent data $(Q,i,a)$ with
the cotangent bundle to the underlying coherent data $(Q,i)$, so
that the moduli stacks of framed $\D$-modules all form pseudotorsors
over the cotangent bundles to certain moduli spaces of complexes of
coherent sheaves on the curve.

A key calculation takes place in Section \ref{TCB section}, where
we identify the space of Koszul data for $\D$-bundles with a
{\em specific} twisted cotangent bundle to the stack of framed torsion
sheaves. We show that the relevant twist comes from the dual
determinant bundle (resulting in Theorem \ref{state main thm}). To
prove this, we use the natural description of the cotangents to quot
schemes to identify the Calogero-Moser moment condition with the
Leibniz rule governing the action map $a$ in the Koszul data. This
calculation can be carried over to describe the moduli space of flat
vector bundles or general holonomic $\D$-modules as twisted
cotangent bundles associated to dual determinant bundles; we hope to
investigate this realization further in future work.

We conclude in Section \ref{explicit} by spelling out very
explicitly the dictionary between our data and the usual CM data on
the affine line.

\subsection{Acknowledgments}
 The authors are
deeply indebted to Victor Ginzburg for generously sharing his ideas.
We would also like to thank Misha Finkelberg, from a lecture of
whose the authors learned of the role of perverse coherent sheaves
in the ADHM classification.

The methods used in the present paper were inspired by two sources.
One of these is the approach of Berest and Chalykh \cite{BC} to
modules over the Weyl algebra (using $A_\infty$-modules over free
algebras). Indeed, one of the goals of this paper is to place their
work in a general dg context and bypass the explicit
use of free algebras, making possible our generalization to
arbitrary (complexes of) $\D$-modules on curves. The authors are
grateful to Yuri Berest for very helpful discussions of \cite{BC}.

A second source for the methods used here is the unpublished work of
Katzarkov-Orlov-Pantev \cite{KOP}.  We are deeply grateful to Tony
Pantev for teaching us how to use Koszul duality and related
methods, and to Katzarkov-Orlov-Pantev for sharing \cite{KOP} with
us.

\section{Noncommutative $\pline$-Bundles and Perverse
Bundles}\label{NC background}
\subsection{Preliminaries on Algebras and Noncommutative $\pline$-Bundles}\label{preliminaries}
Fix a smooth quasiprojective complex curve $X$.
\begin{notation}
Let $\D^1$ denote a bimodule extension of the tangent sheaf $T=T_X$
of $X$ by $\theo_X$ of one of the following two sorts:
\begin{enumerate}
\item $\D^1$ is a Lie algebroid on $X$, i.e. a sheaf of twisted first-order differential operators
$\D^1(L^{\otimes c})$ on some line bundle $L$ on $X$ (where $c\in {\mathbb C}$).
\item $\D^1$ is supported scheme-theoretically on the diagonal $\Delta\subset X\times X$, so it is just given by an extension
of vector bundles on $X$.
\end{enumerate}
\end{notation}
\begin{remark}
Our choice of $\D^1$ above is more restrictive than that required for Theorem \ref{ruled surface thm}.
However, it is easy to see that all the proofs needed for Theorem \ref{ruled surface thm} work with
$\D^1$ replaced by a commutative bimodule $\cE$ as in the statement of the theorem.  The hypotheses
on $\D^1$ chosen above are intended to make statements and proofs easier to follow in the cases of greatest
interest.
\end{remark}
In case (1), we let $\D$ denote the universal enveloping algebroid of $\D^1$, i.e. $\D(L^{\otimes c})$, a sheaf of twisted differential
operators \cite{BB}; in case (2), we let $\D$ denote the quotient
$\on{Sym}(\D^1)/(1-{\mathbf 1})$ of the symmetric algebra of $\D^1$ by the relation that $1\in\on{Sym}^0(\D^1) = \theo$ equals
${\mathbf 1}\in\on{Sym}^1(\D^1) = \D^1$.
This is a filtered $\theo_X$-algebra
(warning: $\theo_X$ is not necessarily a central subalgebra of $\D$!)
 with associated graded algebra
isomorphic to
$\on{Sym}(T_X)$ and with $\D^1$ as the first term in its filtration.
Some  examples: if $\D^1$ is first-order differential operators then $\D$ is the sheaf of
differential operators $\D_X$; if $\D^1 = \D^1(L)$, then $\D$ is the sheaf $\D_X(L)$ of differential
operators
acting on sections of $L$; if $\D^1 = \theo_X\oplus T_X$ (the split extension) then $\D$ is the sheaf $\theo_{T^*X}$ of functions on
the cotangent bundle of $X$.

Let
\bd
\cR = \cR(\D) = \sum_n \D^nt^n \subset \D[t]
\ed
denote the Rees algebra of $\D$; this is a graded algebra with a central
 element
$t\in\cR$ of degree $1$ such that $\cR/t\cR \cong \on{gr}(\D)$.  In particular,
$\cR$ is a graded $\C[t]$-algebra.  The localization
 $\cR[t^{-1}] \cong \D\otimes
{\mbb C}[t,t^{-1}]$ is also a graded ring and we have the
localization functor
$\ell: \cR-\on{mod} \rightarrow \cR[t^{-1}]-\on{mod}$.
Any modules over
graded rings that  we consider will always be assumed to be graded
modules.

We let $\cE$ denote the (formal) microlocalization of $\D$, see e.g.
\cite{AVV}; it is obtained by formally adding to $\D$ power series
in negative powers of vector fields.  For example, if $T_X$ is trivial with nonvanishing section
$\partial$, then $\cE \cong \theo_X(\!(\partial^{\inv})\!)$.
 The ring $\cE$ is also a filtered ring that contains $\D$ as a
filtered subring. Note that $\cE$ is not quasicoherent.
  Given a (left) $\D$-module
$M$, we let $M_\cE =\cE\otimes_\cD M$, which we call (slightly abusively)
the {\em microlocalization of $M$}; similarly, given an
$\cR$-module $N$, we let $N_{\cR(\cE)} = \cR(\cE)\otimes_\cR N$
denote its microlocalization, which
is a graded $\cR(\cE)$-module.
\begin{lemma}\label{localization}
The localization and microlocalization functors $\ell$ and $(-)_{\cR(\cE)}$
are exact; their right adjoints are the forgetful functors back to
$\cR-\on{mod}$.
\end{lemma}
\begin{proof}
For microlocalization this is
Theorem 3.19(2) and Corollary 3.20 of \cite{AVV}.  The statements about
adjoints are standard.
\end{proof}
\begin{remark}\label{microlocal filtration}
As in Section 7 of \cite{solitons}, if $M$ is a filtered $\D$-module
then $M_\cE$ comes equipped with a ``canonical'' filtration.
\end{remark}
\begin{example}
Suppose $\D^1= \theo_X\oplus T_X$ with the central
$\theo_X$-bimodule structure.  Then $\D=\pi_*\theo_{T^*X}$ (where
$\pi: T^*X\rightarrow X$ the projection) is the symmetric algebra of
the tangent sheaf and $\cR$ is the homogeneous coordinate ring of
the ruled surface
$$\PD:=
\proj \cR={\mathbb P}(T^*X\oplus\theo).$$ This is the union of
$\bspec \D=T^*X$ and of a copy of the curve $X=\proj (\cR/t\cR)$
``at infinity". The ring $\cE$ is the ring of formal Laurent series
with poles along the section $X\subset \PD$.
\end{example}
More generally, we will use $\cR$ to construct a noncommutative $\pline$-bundle
\cite{VdB2} as follows.
Let $\on{Gr}\,\cR$ denote the category of graded left $\cR$-modules
that are quasicoherent as $\theo$-modules, and let $\on{gr}\,\cR$
denote its full subcategory of locally finitely generated modules.
We let $\on{Tors}\,\cR$ denote the full subcategory of
$\on{Gr}\;\cR$ whose objects are {\em irrelevant} or {\em locally
bounded modules}: a module $M$ is irrelevant if for every local
section $m$ of $M$ there is some $n_0\in{\mbb Z}$ such that for all
$n\geq n_0$ and all local sections $r$ of $\cR_n$, one has $r\cdot m
= 0$. Geometrically, irrelevant modules are supported
set-theoretically at the irrelevant ideal of $\cR$.
Similarly, $\on{tors}\,\cR$ is the full subcategory of
$\on{gr}\,\cR$ consisting of irrelevant modules.

The subcategories $\on{Tors}\,\cR$ and $\on{tors}\,\cR$ are Serre
subcategories of $\on{Gr}\,\cR$ and $\on{gr}\,\cR$, respectively,
and hence the quotient categories, denoted
$\Qcoh(\PD)=\on{Qgr}\,\cR$ and $\Coh(\PD)=\on{qgr}\,\cR$, exist and
have reasonable properties \cite{AZ1, VdB1}.
\begin{remark}
We will use the notation
$\Qcoh, \Coh$ to emphasize the relationship with derived categories in the
commutative world, since this relationship is absolutely central to our point of view.
\end{remark}
As usual, there is an adjoint pair
of functors \bd (\pi: \on{Gr}\,\cR\rightarrow \Qcoh(\PD),\hskip.2in
\omega: \Qcoh(\PD) \rightarrow \on{Gr}\,\cR) \ed such that the
adjunction $\pi\omega\to Id$ is an isomorphism.   Here $\pi$ is
exact and $\omega$ is left-exact.

 By Theorem 8.8 of \cite{AZ1}
(using the element $t\in\cR$), $\cR$ satisfies the $\chi$-condition;
it follows that Serre's finiteness theorem for cohomology (Theorem
7.4 of \cite{AZ1}) holds for $\cR$.  Moreover, $\cR$ is strongly
noetherian.

All the constructions above sheafify over $X$; in
particular, $\Qgr\,\cR$ forms a sheaf of abelian categories in the
Zariski topology of $X$.

\begin{remark}\label{affine setting}
It is occasionally useful to work algebraically over affine open subsets
of $X$: given an affine open subset $U\subset X$, we may work with the
ring of global sections $\cR(U)$ and carry out the usual constructions for
this ring.  In particular, we will use this at one point later (the proof of Proposition
\ref{Cech facts}).
\end{remark}

\subsection{DG Derived Categories}\label{dg derived cats}
As we mentioned in the introduction, we work throughout not with
derived categories of modules as triangulated categories, but rather
with differential graded (dg) categories which are ``enhancements"
of underlying derived categories. For an excellent overview of dg
categories we refer the reader to \cite{Keller}; detailed treatments
of localization of dg categories can be found also in
\cite{Drinfeld} and \cite{Toen dg}. See also \cite{DAG1} for the
foundations (including a discussion of $t$-structures) of the theory
of stable $\infty$-categories, of which dg categories are the
rational case, and \cite{BLL} for properties of pretriangulated dg
categories (those whose homotopy categories are triangulated).

There are
natural dg versions of essentially all natural constructions with
derived categories, the main technical difference being the need to sometimes
consider dg categories up to quasi-equivalence rather than
equivalence (see the above references as well as \cite{Tamarkin}).
The localization construction for dg categories of Keller, Drinfeld
and To\"en (which is a dg analogue of Dwyer-Kan simplicial
localization) allows one for example to take the dg category of
complexes in an abelian category and localize it with respect to
quasi-isomorphisms, resulting in a {\em canonical} dg category whose
homotopy category is the usual derived category; we refer to
this canonical dg enhancement
as ``the dg derived category."   Another quasi-equivalent dg enhancement
may be obtained from the dg category of injective complexes
(see the discussion in \cite{BLL}).

There are numerous technical advantages in working with dg
enhancements of derived categories rather than with triangulated
categories. Among these we single out three that we use below.
First, the dg enhancements of derived categories satisfy a
homotopical form of descent (or sheaf property) over the \v{C}ech
nerve of a covering (see e.g. Section 7.4 of \cite{BD Hitchin} or
Section 21 of \cite{HirschowitzSimpson}). Second, one can recover
moduli spaces of sheaves from the dg enhanced derived category
\cite{TV}. Finally, dg categories are amenable to explicit
description as modules over the endomorphisms of a compact generator
(see \cite{Keller dct} and the discussion of Theorem \ref{general
framed thm} below). Since we work uniformly throughout the paper
with dg enhancements, we will often abuse terminology and refer
simply to derived categories.

\begin{defn}
Let $D_{dg}(\PD)$ denote the dg enhancement of the bounded derived
category $D^b(\Qgr\,\cR)$ obtained by localizing the dg category of
complexes in $\Qgr\,\cR$ with respect to the class of
 quasi-isomorphisms \cite{Keller,Drinfeld, Toen dg}. We
let $D_{dg,\on{coh}}(\PD)$ denote the full dg subcategory of
$D_{dg}(\PD)$ whose objects have cohomologies in $\on{qgr}\,\cR$.
\end{defn}

\begin{remark}[Calculating in  DG Categories]
The reader who is unfamiliar with
the yoga of dg categories may rest assured that such familiarity is mostly unnecessary for this
paper.  Indeed, although the dg derived category is essential for Theorem \ref{general framed thm}
and for the ``families'' part of Theorem \ref{perverse thm},
the proofs---and, especially, the calculations of objects---depend only on calculations involving
cohomologies and so can be understood in the usual (triangulated) derived category.
\end{remark}

Note that one has an equivalence $\on{Qgr}\,(\on{Sym}T_X) \simeq
\Qcoh(X)$. It follows that one has a base-change functor \bd
i_\infty^*: \Qcoh(\PD)=\on{Qgr}\,\cR \rightarrow \on{Qgr}\,\cR/t\cR
\simeq \Qcoh(X) \ed
and an induced base-change functor $i_\infty^*$ on derived categories.  Both of
these functors may be applied to objects of the abelian category
$\Qgr\,\cR$, and we will write
 $Li_\infty^*$ for the derived functor applied to such objects.
There is also a ``direct image''
functor $(i_\infty)_*$.

\subsection{Perverse Bundles}
We continue to assume that $X$ and $\D$ are as in the previous section.
An object $M$ of $\on{coh}(\PD)$ is said to be {\em zero-dimensional} if
\begin{enumerate}
\item[(a)]
the graded $\cR$-module $\omega M$ has zero-dimensional support on $X$, and
\item[(b)]
the Hilbert function $h_{\omega M}(k) = \on{length}(\omega(M)_k)$ is a
bounded function
of $k$.
\end{enumerate}
  We let $\cT$ denote the full subcategory of $\on{coh}(\PD)$ consisting
of zero-dimensional objects.  We let $\cF$ denote the full subcategory of
$\on{coh}(\PD)$ consisting of objects that have no nonzero zero-dimensional
subobjects.  For every object $N$ of $\on{coh}(\PD)$, there is a unique
maximal zero-dimensional subobject $t(N)$ obtained by forming the sum
$T(\omega N)\subset \omega N$ of all graded $\cR$-submodules of
$\omega N$ that have properties
(a) and (b) above and taking the suboject $t(N) = \pi(T(\omega N))\subset N$.
We thus obtain an exact sequence in $\on{qgr}\,\cR$,
\bd
0\rightarrow t(N)\rightarrow N \rightarrow N/t(N)\rightarrow 0
\ed
with $t(N)\in \cT$ and $N/t(N) \in \cF$.  The pair $(\cT,\cF)$ forms a
torsion pair
in $\on{qgr}\,\cR$, then, and we have:
\begin{prop}[Prop. I.2.1 of \cite{HRS} or Prop. 2.5 of
\cite{Br}]\label{perverse} The full subcategory $\cP$ of
$D_{dg}(\PD)$ given by \bd \cP = \big\{ E\in D_{dg}(\PD) \;\big|\;
H^i(E) =0 \,\text{for $i\notin\{0,1\}$}, H^0(E)\in \cF, H^1(E)\in
\cT\big\} \ed is the heart of a bounded $t$-structure on $D_{dg,
\on{coh}}(\PD)$.
\end{prop}
\noindent Note that by definition a $t$-structure on a
pretriangulated dg category (i.e. one whose homotopy category is
triangulated, \cite{BLL}) is a $t$-structure on its homotopy
category.

We will refer to an object of $\cP$ as a {\em perverse
$\cR$-module}.
An object $N$
 of $\Qcoh(\PD)$ is said to be {\em torsion-free} if it has the
form $\pi(M)$ for some torsion-free module $M$ in $\on{Gr}\,\cR$.
It follows from \cite[S2, p. 252]{AZ1} that in this case $\omega N$ is
a torsion-free graded $\cR$-module.

Let $V$ be a vector bundle on $X$.
\begin{defn}\label{D-bundle def}
A {\em $V$-framed} complex of $\D$-modules is an object $M$ of the
dg derived category $D_{dg}(\PD)$ equipped with an isomorphism
$i_\infty^*(M)\rightarrow V$. A {\em $V$-framed perverse
$\D$-bundle} is a $V$-framed object $M$ of $\cP$ such that
\begin{enumerate}
\item $H^0(M)$ is a torsion-free object of $\Coh(\PD)$, and
\item the $k$th term in the grading $\omega(H^0(M))_k$ is a locally free
$\theo_X$-module of rank $(k+1)\on{rk}(V)$ for all $k\geq -1$.
\end{enumerate}
\end{defn}
\noindent
The second condition is an open condition on perverse $\D$-bundles.  If $\D$ is
commutative, this condition just means that the cohomology sheaf $H^0(M)$ is
isomorphic to $\theo_{\pline}^{\on{rk}(V)}$ on the generic fiber of our $\pline$-bundle
over $X$---in other words, $H^0(M)$  has ``trivial generic splitting type.''
This condition is important when studying $M$ using the Beilinson transform, but
it is also worth noting that it is a good condition that appears ``in nature:" that is, it agrees with
the definition in \cite{solitons} that appears in the study of the KP hierarchy.  Since we will not
need this in the rest of the paper, we will not give the definition of \cite{solitons}---we only
state the following without proof:
\begin{prop}\label{D-bundle equiv}
Suppose that $\D=\D_X(L^{\otimes c})$ is a sheaf of twisted differential operators.  Then the moduli stack of
$V$-framed perverse $\D$-bundles is isomorphic to the moduli stack
of $V$-framed $\D$-bundles in the sense of \cite{solitons}, Definition 3.2.
\end{prop}

We will need the following:
\begin{lemma}\label{perverse basics}
\mbox{}
\begin{enumerate}
\item If $N \in \Coh(\PD)$ is zero-dimensional, then $h_{\omega N}(k)$ is
constant for all $k\gg 0$.
\item Suppose that
$N\in \Coh(\PD)$ is zero-dimensional.  If $i_\infty^*N=0$, then
$Li_\infty^*N=0$.
\item Suppose that $M\in\cP$.  If $H^0(M)$ is torsion-free, then
$ i_\infty^* H^0(M) = Li_\infty^*H^0(M)$.
\end{enumerate}
In particular, if $M$ is a $V$-framed perverse $\D$-bundle, then $
i_\infty^*M = i_\infty^*H^0(M)\cong V$.
\end{lemma}
\begin{proof}
A standard argument using finite generation of $N$ proves (1).

We use the complex
\begin{equation}\label{restriction}
{\mathbf R}: \; \cR(-1)\rightarrow \cR,
\end{equation}
which is quasi-isomorphic to $\cR/t\cR$.  Thus, $ i_\infty^*N \simeq
\pi ({\mathbf R}\otimes \omega N)$. Using (1), it follows that, for
all $k$ sufficiently large,
 multiplication by $t$ from
$\omega N_k$ to $\omega N_{k+1}$ is injective if and only if it
is surjective.  In particular, the homomorphism
$\omega N(-1)\rightarrow \omega N$ has kernel that is a bounded
$\cR$-module if and only if
its cokernel is a bounded $\cR$-module.  This proves (2).

Note that $\omega H^0(M)$ is torsion-free.  It follows that
multiplication by $t$ is injective on $\omega H^0(M)$, and using
\eqref{restriction} to compute $ i_\infty^*H^0(M)$ gives (3).
\end{proof}

\begin{lemma}\label{zero at infinity}
Suppose that $N\in\Coh(\PD)$ is zero-dimensional and $Li^*_\infty N=0$.  Then
(with notation as in Section \ref{preliminaries})
$\ell((\omega N)_{\cR(\cE)}) = 0$.
\end{lemma}
\begin{proof}
If $\D$ is commutative this is well-known.  If $\D$ is a ring of twisted differential operators, then the question reduces to the same
question for the ring of differential operators $\D_X$.  But $\D_X$ has no zero-dimensional modules.
\end{proof}

A perverse $\D$-bundle $M$ has a local numerical invariant called the
{\em (local) second chern class} $c_2(M)$, defined as follows.  Since
$\omega H^0(M)$ is torsion-free, the $V$-framing determines an injective
homomorphism
\bd
\omega H^0(M)/t\omega H^0(M) \hookrightarrow \cR/t\cR\otimes V \cong
\on{Sym}(T_X)\otimes V.
\ed
The cokernel of this homomorphism has finite length $c$ by part (2) of
Definition \ref{D-bundle def}.  Similarly, $\omega H^1(M)_k$ has finite
length
 as an $\theo_X$-module, which is constant for $k\gg 0$; we let
$c'$ denote this length for $k\gg 0$.
\begin{defn}\label{chern class}
With notation as above, $c_2(M) = c+c'$.
\end{defn}

In the commutative case, this definition reproduces the ``local contribution to the second Chern class
of $M$.''  For example, suppose $M$ is given by a complex $\theo_{\PD}\rightarrow Q$, where $Q$ is a
torsion sheaf of length $n$ on the commutative $\pline$-bundle $\PD$.  Then $c_2(M) = n$ regardless of
the map in the complex.  In general, such a ``local $c_2$''---which measures the failure of $H^\bullet(M)$
to be locally free---exists, regardless of whether $X$ is projective (hence our use of the word ``local'').

\section{Resolution of the Diagonal and Beilinson
Transform}\label{Beilinson}
We will next develop the analog, for our noncommutative $\pline$-bundles, of the ``fiberwise
Beilinson transform.''  Although this follows the standard method, it does not seem to appear in the literature in the form
we need.

We also explain the main tool for computing the Beilinson transform, namely,
an analog of \v{C}ech cohomology.  More precisely, we will want to compute direct images from the
$\pline$-bundle $\PD$ to $X$, which we carry out using \v{C}ech cohomology: essentially, this amounts to
covering each $\pline$ by its cover consisting of the affine line (which corresponds to the ring $\D$) and
the ``formal neighborhood of the point at infinity'' (which is captured using the microlocalization $\cE$ of
$\D$).  The \v{C}ech cohomology is easy to calculate in some examples (see Corollary \ref{cech of free}), which will allow
us in Section \ref{identifying moduli} to explicitly compute the Beilinson transform of a perverse bundle.
\subsection{Resolution of the Diagonal}\label{resolution section}
Fix $X$, $\D^1$, and $\cR$ as in Section \ref{preliminaries}.
Let $\Delta$ denote the {\em diagonal bigraded $\cR$-bimodule},
which is given by
$\Delta_{i,j}= \cR^{i+j}$ when $i,j\geq 0$ and $\Delta_{i,j}=0$
otherwise. There is also a larger bigraded bimodule
$\widetilde{\Delta}$ given by $\widetilde{\Delta}_{i,j}= \cR^{i+j}$
for all $i,j$.
\begin{lemma}\label{delta vs delta tilde}
For any object $M$ of $\on{Gr}\,\cR$, the images of the left
$\cR$-modules \bd \pi(\Delta\otimes_\cR M) \;\;\; \text{and}\;\;\;
\pi(\widetilde{\Delta}\otimes_\cR M) \ed are equal in $\Qcoh(\PD)$.
\end{lemma}
\begin{proof}
The quotient module $\widetilde{\Delta}/\Delta$ is
negligible as both a left and a right $\cR$-module.
\end{proof}
Suppose that $M_{\bullet,\bullet}$ is a bigraded sheaf.  We will call $M_{\bullet, \bullet}$ a
{\em bigraded left $\cR$-module} if, for each $k$, $M_{\bullet, k}$ comes equipped with a
structure of graded left $\cR$-module.
\begin{defn}
Given a bigraded left $\cR$-module $M_{\bullet,\bullet}$, we write
$(p_1)_*M = M_{\bullet,0}$, which is a singly-graded left $\cR$-module.
\end{defn}

\begin{lemma}
The functor $\bF:\on{Gr}\,\cR\rightarrow \on{Gr}\,\cR$ defined by
\bd
\bF(M) = (p_1)_*(\wt{\Delta}\otimes_\cR M)= (\wt{\Delta}\otimes_\cR M)_{\bullet,0}
\ed
 is isomorphic
 to the identity
functor.
\end{lemma}
\begin{proof}
This is immediate from the identity
$\widetilde{\Delta}_{i,\bullet}= \cR(i)$ of graded right $\cR$-modules.
\end{proof}

There is an exact sequence of $\cR$-bimodules on $X$:
\begin{equation}\label{alpha sequence}
\cR(-1)\otimes_\theo \D^1\otimes_\theo \cR(-1)
\xrightarrow{\alpha} \cR\otimes_\theo \cR \xrightarrow{m} \Delta
\rightarrow 0.
\end{equation}
Here $m$ is the usual multiplication map and $\alpha = \alpha_1\otimes t -
t\otimes \alpha_2$ is the difference of the two maps given by
\begin{equation}\label{alphas}
\alpha_1: \cR(-1)\otimes \D^1 \rightarrow \cR, \hspace{3em} \alpha_2:
\D^1\otimes \cR(-1)\rightarrow \cR
\end{equation}
 defined by:
\bd
\alpha_1(D_1\otimes Z) = D_1Z, \hspace{5em} \alpha_2(Z\otimes D_2) =  ZD_2.
\ed
  The
kernel of $\alpha$ is easily checked to be the image of
$\cR(-1)\otimes\theo\otimes\cR(-1)$ under the natural inclusion, and so we
obtain:
\begin{lemma}\label{delta res}
The complexes of $\cR$-bimodules
\bd
\delta' = \left[\cR(-1)\otimes\cR(-1)\rightarrow \cR(-1)\otimes \D^1\otimes \cR(-1) \xrightarrow{\alpha}
\cR\otimes\cR\right]
\ed
and
\bd
\delta = \left[\cR(-1)\otimes T_X \otimes \cR(-1) \xrightarrow{\overline{\alpha}}
\cR\otimes\cR\right]
\ed
are quasiisomorphic to $\Delta$.
\end{lemma}
It is important to note that, although $\alpha$ is defined as the
difference of two maps induced from $\theo$-module maps $\alpha_1$
and $\alpha_2$, there is, in general, no such description of
$\overline{\alpha}$.

\subsection{Beilinson Transform}
For any complex of graded $\cR$-modules $\Pbul$, the double complex
$\delta\otimes_{\cR}\Pbul$ is a double complex of
bigraded left $\cR$-modules, and we may form
\begin{equation}\label{beilinson transform}
{\mathbb B}(P_\bullet) = (p_1)_*(\delta\otimes_\cR
\Pbul) \overset{\on{def}}{=} \on{Tot}\left((\delta\otimes_\cR\Pbul)_{\bullet,0}\right);
\end{equation}
here $\on{Tot}$ denotes the total complex.  We call ${\mathbb B}(\Pbul)$ the
{\em Beilinson transform} of $\Pbul$.  Similarly we may replace $\delta$ by $\delta'$ in \eqref{beilinson transform}
and obtain a functor denoted by ${\mathbb B}'$, also called a {\em Beilinson transform}.

\begin{lemma}\label{beilinson exact}
The Beilinson transforms ${\mathbb B}$ and ${\mathbb B}'$ are exact functors.
\end{lemma}
\begin{proof}
Since $\delta$ and $\delta'$  are complexes of projective right $\cR$-modules, ${\mathbb B}$ and ${\mathbb B}'$ are composites of
two exact functors $\delta\otimes_\cR -$ (respectively $\delta'\otimes_{\cR} -$) and $(p_1)_*$.
\end{proof}

\begin{prop}\label{Beilinson of P}
For any complex $\Pbul$ of quasicoherent graded $\cR$-modules, there
are natural morphisms of complexes
\begin{equation}\label{natural morphisms}
{\mathbb B}'(\Pbul)\rightarrow {\mathbb B}(\Pbul)\rightarrow \Pbul.
\end{equation}
If $\Pbul$ is a complex of flat $\cR$-modules, then the cones of \eqref{natural morphisms} have
irrelevant (i.e. locally bounded) cohomologies.
\end{prop}
\begin{proof}
  There are natural maps
\begin{equation}\label{map 1}
\delta'\otimes_\cR \Pbul\rightarrow \delta\otimes_\cR \Pbul\rightarrow \Delta\otimes_\cR \Pbul.
\end{equation}
Suppose that $\Pbul$ is a complex of flat $\cR$-modules; then the maps
\eqref{map 1} are quasiisomorphisms by Lemma \ref{delta res}, as are the natural maps
 \bd
 (p_1)_*(\delta'\otimes_\cR \Pbul) \rightarrow
 (p_1)_*(\delta\otimes_\cR \Pbul) \rightarrow
(p_1)_*(\Delta\otimes_\cR \Pbul).
\ed
Now, the cone of the natural
map
\begin{equation}\label{delta to tilde map}
(p_1)_*(\Delta\otimes \Pbul)\rightarrow
(p_1)_*(\widetilde{\Delta}\otimes \Pbul)
\end{equation}
 is naturally identified with
$(p_1)_*(\widetilde{\Delta}/\Delta \otimes \Pbul)$ which, by Lemma
\ref{delta vs delta tilde}, has irrelevant cohomologies.
\end{proof}
If $M$ is an object of $D_{dg}(\PD)$, the {\em Beilinson transform
of $M$} is defined to be $\pi({\mathbb B}(R\omega(M)))$, respectively
$\pi({\mathbb B}'(R\omega(M)))$
 (here $R\omega$ denotes the total right
derived functor of $\omega$).
\begin{lemma}\label{Beilinson does not depend}
Suppose that $M\rightarrow M'$ is a morphism of complexes in
$\on{Gr}\,\cR$ whose cone has irrelevant cohomologies.  Then the
natural morphisms $\pi{\mathbb B}(M)\rightarrow \pi{\mathbb B}(M')$,
$\pi{\mathbb B}'(M)\rightarrow \pi{\mathbb B}'(M')$
are isomorphisms in $D_{dg}(\PD)$.
\end{lemma}
\begin{proof}
It suffices to check that $\pi{\mathbb B}(N)=0 - \pi{\mathbb B}'(N)$ when $N$ is a module with
irrelevant cohomologies.  Resolving $N$ by locally free $\cR$-modules and
using Lemma \ref{beilinson exact}, this follows from the second
part of Proposition \ref{Beilinson of P}.
\end{proof}

\begin{thm}\label{beilinson thm}
The Beilinson transforms ${\mathbb B}$ and ${\mathbb B}'$ are isomorphic to the identity functor on
 $D_{dg}(\PD)$.
\end{thm}
\begin{proof}
We give the proof for ${\mathbb B}$, the proof for ${\mathbb B}'$ being identical.
By Proposition \ref{Beilinson of P}, there is a natural morphism \bd
\pi((p_1)_*(\delta\otimes R\omega(\Pbul)))\rightarrow \pi\circ
R\omega(\Pbul) = \Pbul. \ed If $\Pbul$ is an object of $\Qcoh(\PD)$
of the form $\pi\cR(k)$, then the cone of the adjunction
$\cR(k)\rightarrow R\omega\circ\pi\cR(k)$ has irrelevant (i.e. locally bounded)
cohomologies \cite[Proposition~7.2]{AZ1}.  Hence, by Lemma \ref{Beilinson does not depend}, we
have an isomorphism $\pi({\mathbb B}(\cR(k)))\rightarrow {\mathbb
B}(\Pbul)$ in $D_{dg}(\PD)$. Now, by Proposition \ref{Beilinson of
P}, we have an isomorphism $\pi{\mathbb B}(\cR(k))\rightarrow
\pi\cR(k) = \Pbul$.  It follows that the natural map ${\mathbb
B}(\Pbul)\rightarrow \Pbul$ is an isomorphism in this case. Since
every object of $D_{dg}(\PD)$ is isomorphic to a complex of direct
sums of objects $\pi\cR(k)$ as above, the
theorem follows by Lemma \ref{beilinson exact}.
\end{proof}

\begin{remark}
Everything in this section and the previous one seems to work
without complications for an arbitrary noncommutative
$\pline$-bundle in the sense of Van den Bergh \cite{VdB2}. This is
interesting to work out in the case of difference operators, which
we plan to carry out in \cite{Toda}.
\end{remark}

\subsection{\v{C}ech Resolution}
In this section, we define an analog of the \v{C}ech resolution for
$\Qcoh(\PD)$.   We will use this to compute Beilinson transforms.

 We begin with some general observations.
For any filtered ring $S$ and its associated Rees ring $\cR(S)$, one has
two functors between the categories of filtered $S$-modules and graded
$\cR(S)$-modules.  The first functor takes a filtered $S$-module $M$ with
filtration $\{M_k\}$ to the graded $\cR(S)$-module $\oplus_k M_k$.
In the other direction, given a graded $\cR(S)$-module $N$, one may invert $t$
and take the degree zero
 term, $f(N) = \ell(N)_0$.   This is naturally a module over
$\cR(S)[t^{-1}]_0 = S$; moreover, $f(N)$ has a filtration given by the
images
of $N_kt^{-k}$ in $f(N)$.  This filtration makes $f(N)$ into a filtered
$S$-module.  These functors are studied in \cite{LvO}.

Suppose now that $M$ is a quasicoherent graded $\cR$-module.  Recall that $\ell(M)$ means the localization of $M$ given by inverting
$t$, and $M_{\cR(\cE)} = \cR(\cE)\otimes_{\cR(\D)} M$, the {\em microlocalization}.
\begin{defn}\label{cech def}
The {\em \v{C}ech complex} of $M$ is the complex
\begin{equation}\label{cech complex}
C(M)= \;\left[ \ell(M) \oplus M_{\cR(\cE)}\xrightarrow{\delta_M}
\ell(M_{\cR(\cE)})\right]
\end{equation}
of $\cR$-modules (implicitly, we compose localization or microlocalization
with the forgetful functor back to $\cR$).  Here the map $\delta_M$ is the difference
of the two natural maps $\ell(M)\rightarrow \ell(M_{\cR(\cE)})$ and
$M_{\cR(\cE)}\rightarrow \ell(M_{\cR(\cE)})$ determined by the adjunctions of Lemma
\ref{localization}.
\end{defn}
\noindent
See the proof of Corollary \ref{cech of free} below for the description of this complex in an example.

As described, we cannot simply apply the functor $\pi$ to get a
\v{C}ech complex $\pi C(M)$ of objects of  $\Qcoh(\PD)$: indeed, the ``microlocal'' terms in $C(M)$ are typically not quasicoherent over $\theo_X$.
To remedy this, we may
replace $C(M)$ with the complex
\begin{equation}\label{Qcoh cech}
\left[\ell(M)\rightarrow \ell(M_{\cR(\cE)})/M_{\cR(\cE)}\right]
\end{equation} which does
consist of objects of $\on{Gr}\,\cR$.  Since the kernel of $M_{\cR(\cE)}\rightarrow \ell(M_{\cR(\cE)})$ is a locally bounded module,
 the images in $\Qcoh(\PD)$ of these complexes are quasi-isomorphic functorially
in $M$.  Hence, we will usually ignore the distinction between (the results of applying $\pi$ to) $C(M)$ and \eqref{Qcoh cech}.
\begin{remark}
As we alluded to above, when $\D$ is commutative, $\ell(M)$ corresponds to the restriction of the quasicoherent sheaf
(corresponding to) $M$ on $\PD$ to $\on{Spec}(\D)$.  Moreover, $M_{\cR(\cE)}$ corresponds to the restriction of $M$ to the
formal neighborhood of the ``curve at infinity" $X_\infty\subset \PD$, and $\ell(M_{\cR(\cE)})$ corresponds to the restriction of $M$
to the ``punctured formal neighborhood'' of $X_\infty$, i.e. the formal neighborhood minus $X_\infty$ itself.  Thus, $C(M)$ is exactly
the \v{C}ech complex corresponding to a cover of the $\pline$-bundle, fiber-by-fiber, by affine lines and formal disks around the points
at infinity in each fiber.
\end{remark}

The \v{C}ech complex comes equipped with a natural morphism $M\rightarrow
C(M)$ given by the adjunction $M\rightarrow \ell(M)$ minus the adjunction
$M\rightarrow M_{\cR(\cE)}$.
\begin{prop}\label{Cech facts}
\mbox{}
\begin{enumerate}
\item The \v{C}ech complex functor is exact.
\item The \v{C}ech complex functor descends to an exact endofunctor of
 $\Qgr\,\cR$.
\item $C(M)$ is a complex of acyclic objects for the functor
$\omega: \Qgr\,\cR\rightarrow\on{Gr}\,\cR$.
\item The natural morphism $M\rightarrow C(M)$ is a quasiisomorphism of
complexes in
$\Qgr\,\cR$.
\end{enumerate}
\end{prop}
\begin{proof}
(1) follows from Lemma \ref{localization}.  (2) follows from Corollaries
4.3.11 and 4.3.12 of \cite{Popescu}, once one observes that inverting $t$
and microlocalization both kill locally bounded modules.

For (3), we first restrict attention to an affine open subset $U$ of $X$
and work as in Remark \ref{affine setting}: that is, we work with the
ring $\cR(U)$ rather than the sheaf of rings $\cR$.  In this case the \v{C}ech
complex as defined in \eqref{cech complex} is a compex in $\Qgr\,\cR(U)$.
We now make
several applications (for the functors  $\ell$, $(-)_{\cR(\cE)}$, and
$\ell(-)_{\cR(\cE)}$) of the following general fact.
\begin{lemma}
Let $\cC, \cD$ be abelian categories and $\cT\subset \cC$ a dense subcategory.
Let $\ell: \cC\rightarrow\cD$ be an exact functor with exact right adjoint
$f:\cD\rightarrow \cC$ so that $\ell(X)=0$ for
all $X\in\on{ob}(\cT)$.  Let
$(\pi:\cC\leftrightarrow \cC/\cT: \omega)$ denote the ``standard'' adjoint
pair of functors.  Let $\overline{\ell}: \cC/\cT\rightarrow \cD$ be the functor
determined by $\overline{\ell}\circ\pi = \ell$ and $\overline{f}$ its right
adjoint.  Then $\overline{f}$ is exact and
$R\omega\circ\overline{f} = \omega\circ\overline{f}$.
\end{lemma}
It follows that \eqref{cech complex} consists of acyclic objects for $\omega$,
and consequently that \eqref{Qcoh cech} also consists of acylic objects.  But
\eqref{Qcoh cech} is now the complex of sections over $U$ of the (sheafified)
\v{C}ech complex, and (3) follows.

To prove (4) it suffices, by (2),
to prove (4) for objects $\pi\cR(k)$, for which see \eqref{coh of free} below.
\end{proof}

\begin{corollary}
If $\Pbul\in D_{dg}(\PD)$, then $C(\Pbul)\simeq R\omega(\Pbul)$.
\end{corollary}

\begin{corollary}[Calculations of \v{C}ech Cohomology]\label{cech of free}
\mbox{}
\begin{enumerate}
\item
We have
\begin{equation}\label{coh of free}
H^\bullet\left((p_1)_*R\omega(\pi\cR(n))\right) = \begin{cases}
\D^n & \text{if $n\geq -1$,}\\
(\cE^{-1}/\cE^n)[-1] & \text{if $n\leq -2$.}
\end{cases}
\end{equation}
In particular, $(p_1)_*R\omega(\pi\cR(-2)) \cong \Omega^1_X[-1]$.
\item For objects $A,B\in D_{dg,\on{qcoh}}(X)$ and $\ell\geq k$, we have
\begin{equation}\label{maps of cplxes}
R\on{Hom}_{\Qgr\,\cR}(\pi\cR(k)\otimes A, \pi\cR(\ell)\otimes B) =
R\on{Hom}_X(A, \D^{\ell-k}\otimes B).
\end{equation}
\end{enumerate}
\end{corollary}
\begin{proof}
Assertion (1)
 follows from a computation of the \v{C}ech cohomology of $\pi\cR(n)$.  Indeed, the \v{C}ech
 complex of $\pi\cR(n)$ is just
 \bd
 \D[t](n) \oplus \cR(\cE)(n) \rightarrow \cE[t](n).
 \ed
 Applying $(p_1)_*$ to this complex (i.e. taking the part in graded degree $0$) gives
 \begin{equation}\label{R coh}
 \D \oplus \cE^n \rightarrow \cE,
 \end{equation}
 where $\cE^n$ is the $n$th term in the filtration of $\cE$.  If $n\geq -1$ then
 \eqref{R coh} is surjective, with kernel $\D\cap \cE^n = \D^n$.  If $n\leq -2$ then
 \eqref{R coh} is injective, with cokernel
 $\cE/(\D + \cE^{n})\cong \cE^{-1}/\cE^n$.

Equation \eqref{maps of cplxes} then follows immediately from
\eqref{coh of free} using the usual
``induction-restriction adjunction.''
\end{proof}

\section{Equivalence and Isomorphism}\label{identifying moduli}
In this section we prove two main consequences of the Beilinson transform of
the previous section.

First, we prove the standard consequence of the Beilinson procedure:
that the derived category of our noncommutative $\pline$-bundle,
$D_{dg}(\PD)$, is quasi-equivalent to a dg category described
completely in terms of ``linear algebra'' or ``quiver data,'' which
we call the dg category $\Kos$ of {\em Koszul data}. Namely, the
Beilinson transform shows that $D_{dg}(\PD)$ is generated by the
analog of the exceptional collection $\{\Oo(-1),\Oo\}$ and hence may
be described via modules over the endomorphism algebra of this
collection.
  In our setting, when everything is done in families over the curve $X$, this
amounts to describing an object of $D_{dg}(\PD)$ in terms of a pair of quasicoherent complexes
on the curve $X$ together with the data encoding the lift to the $\pline$-bundle, i.e.
an action of $\D^1$ from one complex on $X$ to the other.

Second, we describe explicitly (Theorem \ref{perverse thm}) which Koszul data (objects of $\Kos$)
 correspond to perverse bundles in $D_{dg}(\PD)$.  This amounts to doing a \v{C}ech calculation starting
 from a perverse bundle, and, conversely, computing the cohomologies of an object of $D_{dg}(\PD)$ starting
 from an object of $\Kos$.  As we explain below, once one has carried through these calculations at the level of
 objects, the framework of dg categories gives the equivalence of the moduli essentially ``for free.''

All tensor products in this section are taken over $\theo_X$.

\subsection{DG Category of Koszul Data}\label{defn of Kos}
Let $P$ denote the sheaf of rings on $X$ consisting of $2\times 2$ upper triangular matrices whose diagonal coefficients lie in
$\theo_X$ and whose upper-left entry lies in $\D^1_X$:
\bd
P = \left\{
\begin{pmatrix}
f_1 & D\\
0 & f_2
\end{pmatrix}
\; \Bigg| \;\; f_1, f_2 \in \theo_X, D\in \D^1 \right\}. \ed This
algebra is the version, for our noncommutative $\pline$ bundle
$\PD$, of the endomorphism (or self-Ext) algebra of the Beilinson
generator $\Oo(-1)\oplus \Oo$ on $\pline$ (see below).

We let $\on{Kos}$ denote the abelian category consisting of quasicoherent sheaves $M$ together with a map
$P\otimes_{\theo_X} M\rightarrow M$ that defines a $P$-module structure (and whose morphisms are $P$-linear maps).
It is immediate that specifying such a module is the
same as giving a pair of quasicoherent $\theo_X$-modules $M_{-1}, M_0$  (our choice of labelling will seem more natural in light of
Theorem \ref{general framed thm}) together with an $\theo_X$-linear map
$\D^1\otimes_{\theo_X} M_{-1} \xrightarrow{a} M_0$, an {\em action map}.
\begin{defn}
The dg category of {\em Koszul data} is the dg derived category $\Kos$ of the abelian category $\on{Kos}$.
\end{defn}
Objects of $\Kos$ are thus given by pairs $(M_{-1}, M_0)$ of quasicoherent complexes on $X$ together with an
action map $a: \D^1\otimes M_{-1}\rightarrow M_0$.
Let $V\in D(X)$ be a coherent complex on $X$.
An object $(M_{-1}, M_0, a)$ of $\Kos$ is {\em $V$-framed} if it comes equipped
with a quasi-isomorphism $\on{Cone}(M_{-1}\xrightarrow{a|_{M_{-1}}} M_0)
\simeq V$; here $a|_{M_{-1}}$ means the restriction of $a$ to $M_{-1}\subset \D^1\otimes M_{-1}$.

\subsection{General Equivalence Theorem}
Let $G  = \pi \cR \oplus \pi\cR(1)$ in $D_{dg}(\PD)$ (of course, this is the image of an object of $\Coh(\PD)$).
It follows from Corollary \ref{cech of free} that
the sheaf of dg algebras $\Hom_{D_{dg}(\PD)}^\bullet(G,G)$ is canonically quasi-isomorphic to the sheaf of algebras $P$
defined above (that is, concentrated in cohomological degree $0$).

Given an object $\cM$ of $D_{dg}(\PD)$, we get a complex of (left!) $P$-modules $M$ by setting
$M = \Hom^\bullet_{D_{dg}(\PD)}(G, \cM)$.
It is instructive to apply this procedure to an output of the Beilinson functor ${\mathbb B}$, i.e. to the total complex associated to
a double complex
\begin{equation}\label{dbl cplx}
\cM: \;\; \pi\cR(-1)\otimes T_X\otimes M_{-1}\rightarrow \pi\cR \otimes M_0
\end{equation}
(that is, $M_{-1}$ and $M_0$ are quasicoherent complexes on $X$).
It follows from Corollary \ref{cech of free} that:
\begin{lemma}\label{kos of bei}
For $\cM$ as in \eqref{dbl cplx},
\bd
\Hom^\bullet(\pi\cR(1), \cM)\simeq M_{-1} \;\; \text{and} \;\;
\Hom^\bullet(\pi\cR, \cM)\simeq M_0.
\ed
\end{lemma}

We now have:
\begin{thm}\label{general framed thm}
\mbox{}
\begin{enumerate}
\item
There is a quasi-equivalence of differential graded categories
from $D_{dg}(\PD)$  to $\Kos$ that takes an object $\cM$ of
$D_{dg}(\PD)$ to $\Hom^{\bullet}_{D_{dg}(\PD)}(G, \cM)$.
\item
 Under this quasi-equivalence, if an object $\cM$
of $D_{dg}(\PD)$ is identified with a triple $M = (M_{-1}, M_0, a)$ of $\Kos$, then a
choice of $V$-framing of $\cM$ corresponds to a choice of $V$-framing
of $M$.
\end{enumerate}
\end{thm}
\begin{proof}
By definition of quasi-equivalence \cite{Keller}, it suffices to
define a functor of dg categories that induces an equivalence of the
homotopy categories (as triangulated categories).
   The functor described above at the level of objects defines a functor from $D_{dg}(\PD)$ to $\Kos$ by a standard procedure
(see Section 4 of \cite{Keller ddgc} or Section 8.7 of \cite{Keller
dct}). Moreover, by Theorem 4.3 of \cite{Keller ddgc} (see also the
theorem in Section 8.7 of \cite{Keller dct}), if $X$ is affine---so
that the category $\Kos$ is a (derived) module category for a ring,
not just a sheaf of rings---then, to conclude that this functor
induces an equivalence of the homotopy categories, it suffices to
check that the compact object $G$ is a generator of $D_{dg}(\PD)$.
But this is immediate from Theorem \ref{beilinson thm}.  Finally,
for projective $X$, we obtain a compatible collection of such
functors for all open subsets of $X$. By the above discussion, these
functors are quasi-equivalences for all affine open subsets of $X$.
It then follows from effective descent for dg categories that the
functor $\Hom^\bullet(G, -)$ gives a quasi-equivalence over $X$ as
well. More precisely, the dg categories $D_{dg}(\PD)$ and $\Kos$ are
homotopy limits of the categories over affine open subsets of $X$ --
they are obtained by totalizing the cosimplicial dg categories
associated to the \v{C}ech nerve of an affine cover. (See \cite{BD
Hitchin} (Section 7.4) or \cite{HirschowitzSimpson} (Section 21) for
descent of dg categories of sheaves.)

We now prove the second part of the theorem.
We first note the following.  Consider the complex of bigraded $\cR$-bimodules $\delta$ in Section \ref{resolution section}.
Take the associated graded complex with respect to the left-hand grading.  The description of the map $\overline{\alpha}$ in terms of
$\alpha_1$ and $\alpha_2$ shows that, after passing to associated gradeds with respect to the left-hand grading,
$\overline{\alpha}$ is identified with $\alpha_1\otimes t$.

Consider now a double complex of the form \eqref{dbl cplx}, an output of the functor
${\mathbb B}$; by Lemma \ref{kos of bei},
the corresponding Koszul data are $(M_{-1}, M_0, a)$.
Applying $i_\infty^*$ to
\eqref{dbl cplx} gives
\bd
\left[\pi(\on{Sym}T_X)(-1) \otimes T_X \otimes  M_{-1} \longrightarrow
\pi\on{Sym}(T_X)\otimes M_0\right].
\ed
By the conclusion of the previous paragraph, this map is just the tensor product of the identity map $\pi\on{Sym}\,T_X(-1)\otimes T_X \rightarrow \on{Sym}\,T_X$ with $a|_{M_{-1}}$.
Under the identification of $\Qcoh(X)$ with $\Qgr\,\on{Sym}(T_X)$ (more precisely, of their derived categories), this is then identified with $a|_{M_{-1}}$.  Part (2) of the theorem
is then immediate from the two definitions of $V$-framing.
\end{proof}

\subsection{Calculating in $\Kos$}
Let $V$ be a sheaf on $X$.  We want to calculate some $V$-framed objects of $\Kos$ that have a particularly simple form:
their components $M_{-1}$ are (quasi-isomorphic to) sheaves in cohomological degree $1$.

Suppose $M$ is a $V$-framed object of $\Kos$.  As we have remarked above, since our sheaf of rings $P$ has two projections onto
$\theo_X$, $M$ gives two quasicoherent complexes $M_{-1}$ and $M_0$; more precisely, we get two functors to
$D_{dg}(X)$.  Suppose that $M_{-1}$ is quasi-isomorphic to a sheaf in cohomological degree $1$.  Letting
$Q=H^1(M_{-1})$, we then get a canonical quasi-isomorphism $M_{-1}\rightarrow Q[-1]$.   Furthermore, the
$V$-framing of $M$---that is, the choice of quasi-isomorphism
$\on{Cone}(M_{-1}\rightarrow M_0)\simeq V$---yields a choice of quasi-isomorphism $M_0\rightarrow \on{Cone}(V\rightarrow Q)$ (with $V$ in cohomological degree $0$),
i.e. up to quasi-isomorphism $M_0$ corresponds to a choice of morphism $V\xrightarrow{i} Q$.  The action map
$a:\D^1\otimes M_{-1}\rightarrow M_0$ then yields a unique map
$a: \D^1\otimes Q[-1]\rightarrow \on{Cone}(V\rightarrow Q)$ in $D_{dg}(X)$.  This map has the property that the composite
\begin{equation}\label{composite}
Q[-1] \rightarrow \D^1\otimes Q[-1] \rightarrow \on{Cone}(V\rightarrow Q)
\end{equation}
equals the ``natural'' map $Q[-1]\rightarrow \on{Cone}(V\rightarrow Q)$.

Suppose, for the moment, that $Q$ is a quasicoherent sheaf on $X$, $V\xrightarrow{i} Q$ is a choice of homomorphism,
and $a:\D^1\otimes Q[-1]\rightarrow \on{Cone}(V\xrightarrow{i} Q)$ is a choice of map in $D_{dg}(X)$ whose composite
\eqref{composite} equals the natural map $Q[-1]\rightarrow \on{Cone}(V\xrightarrow{i} Q)$.
  We will refer to $(Q,i,a)$ as a
{\em Koszul triple} on $X$.
 We will define the set of {\em morphisms} of triples $(Q,i, a)\rightarrow (Q', i', a')$ to consist of morphisms
$Q\xrightarrow{\phi} Q'$ so that the diagram
\bd
\xymatrix{
\D^1\otimes Q[-1]\ar[r]^{a} \ar[d]^{1_{\D^1}\otimes\phi} & \on{Cone}(V\xrightarrow{i} Q)\ar[d]^{1_V\times \phi}\\
\D^1\otimes Q'[-1]\ar[r]^{a'} & \on{Cone}(V\xrightarrow{i'} Q)
}
\ed
commutes in $D_{dg}(X)$.
 This defines a category
of triples $(Q,i,a)$.  Similarly, we let $\Kos^0(V)$ denote the category of $V$-framed objects $M$ in $\Kos$ such that
$M_{-1}$ is quasi-isomorphic to a sheaf in cohomological degree $1$;  Homs are
morphisms in $\Kos$ that are compatible with the $V$-framing.
\begin{prop}\label{special framed remark}
There is an equivalence of categories between:
\begin{enumerate}
\item
The category $\Kos^0(V)$.
\item The category of triples $(Q,i,a)$.
\end{enumerate}
\end{prop}
\begin{proof}
The above constructions define a functor from $\Kos^0(V)$ to the category of Koszul triples.  It is faithful, since our analysis of
$M_{-1}$ and $M_0$ shows that any map of Koszul data $M\rightarrow M'$ is determined by the induced map
$Q = H^1(M_{-1})\rightarrow H^1(M'_{-1}) = Q'$.

To see that this functor is full and essentially surjective, we do the following.  Given $(Q,i,a)$, choose representatives
$M_{-1}$ and $M_0$ of $Q[-1]$ and $V\rightarrow Q$ in $D_{dg}(X)$, respectively, for which a lift
$a:\D^1\otimes M_{-1} \rightarrow M_0$ exists (for example, choose a complex of injectives representing $M_0$).
Applying the functor to the Koszul data $(M_{-1}, M_0, a)$ returns $(Q,i,a)$---that is,
the functor is essentially surjective.
Furthermore, given another triple $(Q', i', a')$, we may choose complexes of injectives $M_{-1}'$ and $M_0'$ quasi-isomorphic to
$Q'$ and $V\rightarrow Q'$, so that $a'$ lifts to an action map $\D^1\otimes M_{-1}'\xrightarrow{a'} M_0'$ and $(M_{-1}', M_0', a')$ is
mapped to $(Q', i', a')$.
Then any map $Q\rightarrow Q'$ lifts to a map $M_{-1}\xrightarrow{\phi} M_{-1}'$.  It then
suffices to choose a map $M_0\rightarrow M_0'$ compatible with $\phi$ and the given quasi-isomorphisms to $V$, but the existence
of such a map is guaranteed by injectivity of $M_0'$ using an exact sequence argument.
\end{proof}

\noindent
We will use this description in
Theorem \ref{perverse thm} below.

\subsection{Moduli of Perverse Bundles}
Fix a vector bundle $V$ on $X$.
\begin{defn}\label{stack defs}
Let $\PS(X,V)$ denote the moduli stack of pairs $(Q,i)$ consisting of a coherent torsion sheaf
$Q$ on $X$ and a homomorphism $i:V\rightarrow Q$; we call the component parametrizing pairs $(Q,i)$ with $Q$ of
length $n$  the {\em $n$th perverse symmetric power} of $X$.

 In the case $\D^1 = \theo_X\oplus T_X$ (commutative), we let $\PB(X,V)$ denote the moduli stack of $V$-framed
perverse $\D$-bundles.  In the case when $\D^1 = \D^1_X$ consists of first-order differential operators, we let
$\CM(X,V)$ denote the moduli stack of $V$-framed perverse $\D$-bundles.
\end{defn}

We will now describe the moduli stacks $\PB(X,V)$ and $\CM(X,V)$ in terms of stacks that ``live over'' $\PS(X,V)$.

\begin{thm}\label{perverse thm}
The moduli stacks for the following data are isomorphic:
\begin{enumerate}
\item $V$-framed perverse $\D$-bundles $\cM$ with second Chern class $c_2$.
\item triples $(Q,i,a)$ where $Q$ is a torsion sheaf on $X$ of
length $c_2$, $i:V\rightarrow Q$ is a homomorphism of
$\theo$-modules, and \bd a: \D^1\otimes Q[-1] \rightarrow
\on{Cone}(V\xrightarrow{i} Q) \ed
 is a morphism
in the dg derived category of $\theo$-modules such that the induced
morphism \bd \theo\otimes Q[-1]\rightarrow \on{Cone}(V\xrightarrow{i} Q) \ed
is the natural one.
\end{enumerate}
\end{thm}
\begin{remark}
It is instructive to consider what the theorem says in case $\cR$ is
commutative and $\cM$ is a rank $1$ torsion-free sheaf. Let $X_\infty$
denote the ``divisor at infinity'' in $S= \on{Proj}(\cR)$ and
$i_\infty: X_\infty \rightarrow S$ the inclusion.  Let
$p:S\rightarrow X$ denote the projection.  Note that $\on{Qgr}(\cR)
\simeq \Qcoh(S)$.

By condition (2) of Definition \ref{D-bundle def}, $\cM$ comes equipped with
an injective homomorphism $\cM\hookrightarrow p^*V$.  We let $\cQ$ denote the
cokernel; this gives an exact sequence
\begin{equation}\label{rank 1}
0\rightarrow \cM\rightarrow p^*V \rightarrow \cQ\rightarrow 0.
\end{equation}
We have
 $M_0 = {\mathbb R}p_*(\cM)$ and $M_{-1}= {\mathbb R}p_*(\cM(-X_\infty))$.  We also
set $Q=p_*\cQ$; this is a torsion sheaf of length $c_2(\cM)$.
  Using the sequence \eqref{rank 1}, one finds that
$M_0\simeq [V\xrightarrow{i} p_*\cQ]$ and $M_{-1}\simeq
p_*\cQ(-X_\infty)\cong Q$.
\end{remark}

\begin{proof}[Proof of Theorem \ref{perverse thm}]
In light of Proposition \ref{special framed remark}, it
suffices to check that, under the
correspondence of Theorem \ref{general framed
thm} above, the perverse $\D$-bundles correspond to Koszul data
$(M_{-1}, M_0, a)$ in which $M_{-1}= (p_1)_*R\omega(M(-1))$ is quasi-isomorphic to  $Q[-1]$
for $Q$ a torsion sheaf of length
$c_2$.

The strategy of the proof, then, is simple: we start with a perverse bundle $\cM$ and compute
\bd
\Hom^\bullet(G,\cM)  = M_{-1}\oplus M_0 \simeq (p_1)_*R\omega(\cM) \oplus (p_1)_*R\omega(\cM(-1))
\ed
and show that it has the desired form.  Then, going the other way, we start with a triple $(Q,i,a)$---equivalently by
Proposition \ref{special framed remark}, with framed Koszul data---and compute the corresponding object of
$D_{dg}(\PD)$.

{\em Correspondence for ${\mathbb C}$-points.}\;
We first give the proof at the level of ${\mbb C}$-valued points.  Later we will explain that this is already
enough to give the equivalence in
families, i.e. an isomorphism of stacks.

Suppose first that $\cM$ is a $V$-framed
perverse $\D$-bundle.  In order to check that the corresponding triple $(M_{-1}, M_0, a)$ of
Koszul data has the desired property, it suffices to compute the ``derived direct image''
$(p_1)_*R\omega \cM(-1)$ and prove that it is a torsion sheaf (in cohomological degree $1$) of length $c_2$.
To carry  out this calculation,
we use the exact triangle
\begin{equation}\label{comm M exact}
H^1(\cM)(-1)[-1] \rightarrow \cM(-1) \rightarrow H^0(\cM)(-1)\xrightarrow{[+1]}
\end{equation}
in the triangulated derived category.
  By Lemma \ref{perverse basics}, the natural map
\bd
V\simeq i_\infty^*\cM(-1)\rightarrow Li_\infty^*H^0(\cM)(-1)
\ed
 is a quasi-isomorphism, and
it follows that $Li_\infty^* H^1(\cM)(-1) = 0$: in other words, ``the zero-dimensional object
$H^1(\cM)(-1)$ is supported away from $X_\infty$.''

Write $H = \omega H^1(\cM)(-1)$.  By Lemma \ref{zero at infinity} and the conclusion of the previous paragraph,
 $\ell(H_{\cR(\cE)})=0$.
 Thus, using Formula \eqref{Qcoh cech} for the \v{C}ech complex $C(H)$, $C(H)$  is given by $C(H) \simeq
H[t^{-1}]$ and $(p_1)_*C(H)$ is a torsion coherent sheaf on $X$ of  length $c'$ (notation as in
Definition \ref{chern class}).
For $L = H^0(\cM)(-1)$, the argument of
 Proposition 5.8 of \cite{solitons} (more specifically, the exact sequence
(5.2)) applies here to prove that $H^0((p_1)_*C(L))=0$ and $H^1((p_1)_*C(L))$
 is a coherent $\theo_X$-torsion
sheaf of length  $c$ (notation as in Definition
\ref{chern class}).

Taking account of the cohomological degrees of the ``direct images'' of the
 left- and right-hand terms in \eqref{comm M exact},
 it is now immediate from
 that $(p_1)_*R\omega(\cM(-1))[1]$ is concentrated in a single cohomological degree. It follows that
 $(p_1)_*R\omega(\cM(-1))[1]$ is a torsion $\theo_X$-module of
length $c+c' = c_2$, as desired.  Summarizing, we have proven that a $V$-framed perverse $\D$-bundle $\cM$ gives a triple of
Koszul data as claimed in the statement of the theorem.

We will now prove the converse: Koszul data that reduce, via Proposition \ref{special framed remark}, to a triple $(Q,i,a)$  as in part (2) of the theorem correspond, under Theorem \ref{general framed thm}, to
a perverse bundle with the prescribed $c_2$.

 Suppose we are given Koszul data $(M_{-1}, M_0, a)$ that correspond to a triple
$(Q,i,a)$ as in part (2) of the present theorem.  We set
\begin{equation}\label{dg from curve to surface}
\cM  =\on{Tot}\big[\pi\cR(-1)\otimes T_X\otimes M_{-1}
\longrightarrow \pi\cR\otimes M_0\big],
\end{equation}
the object of $D_{dg}(\PD)$ corresponding to $(M_{-1}, M_0, a)$
(see Lemma \ref{kos of bei}, and note that $\pi\cR\otimes M_0$ lies in cohomological degree $0$, not $1$).   To compute whether $\cM$ is a perverse bundle, it suffices to compute in the triangulated
category---that is, we may compute cohomologies of $\cM$ using the quasi-isomorphic object
\begin{equation}\label{from curve to surface}
\overline{\cM}  =\on{Tot}\big[\pi\cR(-1)\otimes T_X\otimes Q[-1]
\longrightarrow \pi\cR\otimes \on{Cone}(V\rightarrow Q)\big]
\end{equation}
in the triangulated derived category (again, note the normalization, that the right-hand term $\pi\cR\otimes\on{Cone}(V\rightarrow Q)$ of the double complex lies in ``horizontal cohomological degree $0$'').   In particular, one sees from this description of
$\overline{\cM}$ that $H^k(\cM)=0$ for $k\neq 0, 1$.
Furthermore, over
$U=X\setminus\on{supp}(Q)$, \eqref{from curve to surface} immediately reduces to
\begin{equation}\label{generic}
\overline{\cM}|_U =
V\otimes\pi\cR|_U;
\end{equation} in particular, this implies the rank condition of part (2) of
Definition \ref{D-bundle def}---the ``locally free'' part of the condition is immediate provided $H^0(\cM)$ is a torsion-free
object of $\Qcoh(\PD)$, which we will prove below.

We next prove that $H^1(\cM)$ is zero-dimensional.
By \eqref{generic}, $\omega H^1(\cM)$ is supported over $\on{supp}(Q)\subset X$.
Moreover, applying $i_\infty^*$ to  \eqref{from curve to surface}, we find that $i_\infty^*H^1(\cM)$ is the cokernel of a map
\begin{equation}\label{restriction2}
\pi(\on{Sym}\, T_X \otimes V)\oplus \big((\on{Sym}\, T_X(-1))\otimes T_X \otimes Q\big) \longrightarrow \pi \on{Sym}\, T_X \otimes Q.
\end{equation}
By an argument similar to the proof of part (2) of Theorem \ref{general framed thm}, this map is surjective; in particular,
$i_\infty^* H^1(\cM)=0$.  It follows that $\omega H^1(\cM)$ is $\theo$-coherent.  Since $H^1(\cM)$ is also $\theo_X$-torsion (supported on
$\on{supp}(Q)$), it is zero-dimensional.

By the previous paragraph and Lemma \ref{perverse basics}, $Li_\infty^*H^1(\cM) = 0$, and hence, since
\eqref{restriction2} comes equipped with a quasi-isomorphism to $\on{Sym}\, T_X \otimes V$, we get a quasi-isomorphism of $i_\infty^*\cM$ with $V$.  Hence
$\cM$ is $V$-framed.

To see that $\cM$ is a perverse bundle, it remains only to prove that $H^0(\cM)$ is torsion-free.
The description \eqref{from curve to surface}  gives us an exact triangle
\bd
\pi\cR\otimes \on{Cone}(V\rightarrow Q) \longrightarrow \overline{\cM} \longrightarrow \pi\cR(-1)\otimes T_X \otimes Q\xrightarrow{[1]}
\ed
(the lack of shifts may look unexpected, but results from the shifts implicit in our description of $\tilde{M}$ via \eqref{from curve to
surface}).
The associated long exact sequence
takes the form
\bd
0\rightarrow H^0\big((\pi\cR\otimes\on{Cone}(V\rightarrow Q)\big) \rightarrow H^0(\cM)
\rightarrow H^0(\pi\cR(-1)\otimes T_X \otimes Q).
\ed
Applying $\omega$, this remains left exact.
Moreover, the left-hand term $\omega\,H^0\big(\pi\cR\otimes\on{Cone}(V\rightarrow Q)\big)$ is a subobject of $\cR\otimes V$, hence
is torsion-free.
Consequently,
if $\omega\,H^0(\cM)$ is not torsion-free, then its torsion submodule $\tau$ maps
injectively to a submodule of
\bd
\omega\,H^0(\pi\cR(-1)\otimes T_X\otimes Q)=\cR(-1)\otimes T_X\otimes Q;
\ed
 in particular,
the torsion submodule is a direct summand of $\omega\,H^0(\cM)$ and is isomorphic to
a submodule of $\cR(-1)\otimes T_X\otimes Q$.  Now $i_\infty^*H^0(\cM)$ is
torsion-free,
so $i_\infty^*\pi(\tau)=0$.  It follows that $t \cdot \tau_k =\tau_{k+1}$ for $k\gg 0$.
But it is clear that no nonzero submodule of $\cR\otimes T_X\otimes Q$ has this
property.  This proves that $H^0(\cM)$ is torsion-free, thus completing the proof that $\cM$ is a $V$-framed
perverse $\D$-bundle.

We thus have the desired equivalence at the level of ${\mathbb C}$-valued points of the stack.

\vspace{.7em}

\noindent
{\em Correspondence in families.}\;  This is essentially a formal consequence
of what we have already proven.  More precisely, following
 \cite[Definition~2.1.8 and Proposition~2.1.9]{Lieblich} (or, more generally, \cite{TV}) there
is an intrinsic notion of $S$-object in a dg category (where $S$ is
a locally noetherian scheme).  A quasi-equivalence of dg categories induces an
equivalence of the associated categories fibered in groupoids over schemes.
The only additional information that one would like is that, if
one has a ``classical'' notion of flat family, the two notions coincide.

For families
${\mathcal M}$
of $V$-framed perverse $\D$-bundles for which $H^0({\mathcal M}_s)=0$
for every geometric point $s\in S$ this is straightforward
(note, for example, that this applies to {\em every} family when
$\D$ is a TDO algebra).  Such a family certainly satisfies the condition
$\on{Ext}^i({\mathcal M}_s,{\mathcal M}_s) = 0$ for $s\in S$ and
$i<0$ that appears
in the definition \cite{Lieblich} of a universally gluable
family.  Conversely, suppose
 ${\mathcal M}$ is a universally gluable $S$-object of
$D_{dg}(\PD)$ all of whose geometric fibers ${\mathcal M}_s$ are
($V$-framed) perverse $\D$-bundles and such that $H^i({\mathcal
M}_s) = 0$ for $i\neq 0$ (note that here ${\mathcal M}_s$ means
derived restriction). It follows by \cite[Lemma~4.3]{bridgeland
fourier} that $H^i({\mathcal M})=0$ for $i\neq 0$ and that
$H^0({\mathcal M})$ is $S$-flat.
\end{proof}

\section{Calogero-Moser Spaces as Twisted Cotangent Bundles}\label{TCB section}

In this section we identify precisely the moduli stacks of perverse
bundles $\PB$ on $T^*X$ and $\D$-bundles $\CM$ and the relation
between them. Recall (Definition \ref{stack defs}) that ${\PS}(X,V)$ denotes the moduli stack of
pairs $(Q,i)$ consisting of a coherent torsion sheaf $Q$ on $X$ and
a homomorphism $i:V\rightarrow Q$. For fixed $(Q,i)$ we may now
consider its possible extensions to data  $(Q,i,a)$ which by Theorem
\ref{perverse thm} describe either perverse bundles or $\D$-bundles.
In the first case, $a$ is completely determined by a map
$\ol{a}:T_X\ot Q[-1]\to\on{Cone}(V\xrightarrow{i} Q)$, while in the latter
case $a$ is a map $\D^1\ot Q[-1]\to\on{Cone}(V\xrightarrow{i} Q)$ with fixed
restriction to $\theo\ot Q[-1]$. It is clear that the latter data
form a pseudo-torsor over the former: in other words, there is a
simply transitive action of the perverse bundle structures on
$(Q,i)$ on the (possibly empty) set of $\D$-bundle structures on the
same $(Q,i)$ (both over $\C$ and in families).

It is easy to compute from the deformation theory of pairs $(Q,i)$
that the data $\ol{a}$ are dual to the tangent space to $\PS(X,V)$
at $(Q,i)$. However we would like to be more precise in the
identification of the moduli of $\D$-bundles as (the stack analog
of) a twisted cotangent bundle in the sense of \cite{BB}, i.e. as
(pseudo)torsor for the cotangent bundle equipped with a compatible
symplectic form. Recall that two standard constructions of a twisted
cotangent bundle are as the affine bundle of connections on a line
bundle \cite{BB}, and as a ``magnetic cotangent bundle" \cite{RSTS},
a Hamiltonian reduction at a one-point coadjoint orbit. The two
constructions agree in the case of reduction at a coadjoint orbit
coming from a character of the group, which produces the affine
bundle of connections on the line bundle obtained by descending the
trivial line bundle with equivariant structure given by the
character.

We will describe both the perverse bundle and the $\D$-bundle spaces
as Hamiltonian reductions at (zero and nonzero, respectively)
integral one-point orbits from the cotangent bundle of the relevant
Quot-type scheme. Recall from \cite{BD Hitchin} (Section 1.2) that
the very definition of the cotangent bundle of a stack is naturally
given as a Hamiltonian reduction (at $0$) of the cotangent bundle of
a cover. Thus perverse bundles naturally form the cotangent bundle
of a Quot scheme, while $\D$-bundles form the space of connections
of a line bundle which we identify with the dual of the determinant
line bundle.

\begin{remark}\label{dg cotangent stack}
In addition to its concrete flavor, the description of $\PB$ and
$\CM$ via Hamiltonian reduction from a Quot scheme also has an
 implicit advantage. Namely,
 it is easily adapted to capture the cotangent and
twisted cotangent stacks $\PB$ and $\CM$ as derived (or dg) stacks,
rather than naively as stacks.
\end{remark}

Let $\wt{\PS}_n(X,V)$ denote the moduli space of triples
$(Q,i,\ol{q})$ where $(Q,i)\in \PS_n(X,V)$ and $\ol{q}:\C^n\simeq
H^0(Q)$ is a basis of sections of $Q$.
The datum of $\ol{q}$ is
equivalent to that of a surjection $q:\C^n\ot\Oo\to Q$ (that induces
an isomorphism $\C^n\rightarrow H^0(Q)$), so that $\wt{\PS}_n(X,V)$
is the fiber product (over length $n$ sheaves $Q$) of $\PS_n(X,V)$
with an open set of the
 Quot scheme parametrizing
$n$-dimensional quotients $q$.  As an immediate consequence of this description one has:
\begin{lemma}
The space $\wt{\PS}_n(X,V)$ is a smooth variety.
\end{lemma}
 Note that
$$\PS_n(X,V)=\wt{\PS}_n(X,V)/GL_n,$$
where $GL_n$ acts by changing the basis of
$H^0(Q)$.

\begin{lemma} The cotangent bundle of $\wt{\PS}_n(X,V)$ is the
stack of quadruples $(Q,i,q,r)$ with $(Q,i,q)$ as above and
\begin{equation}\label{r}
r: Q[-1]\to \on{Cone}(V\oplus \Oo^n \xrightarrow{(i,q)} Q)\ot
\Omega.
\end{equation}
\end{lemma}
\begin{proof}
First suppose that $X$ is projective.  Then a standard calculation
gives \bd T_{(Q,i)}\wt{\PS}_n(X,V) =
\on{Hom}(V\oplus\theo^n\xrightarrow{(i,q)} Q, Q) \ed and the lemma
follows by Serre duality.  If $X$ is only quasiprojective, we
complete $X$ to a projective curve $Y$ (and extend $V$) and apply
the above argument to $Y$; the lemma then follows by noting that all
of our data $(Q,i,q,r)$ are in fact supported on $X\subseteq Y$.
\end{proof}
\begin{remark}
From this point forward, we may always assume that $X$ is
projective.  Indeed, if not, then we may complete $X$ to a
projective curve $Y$ and extend $V$ to a vector bundle
$\overline{V}$ on $Y$.  Since $\PS(X,V)$ and $\wt{\PS}(X,V)$ are
open substacks of $\PS(Y,\overline{V})$ and
$\wt{\PS}(Y,\overline{V})$ respectively, the conclusions of
Proposition \ref{PB prop} and Theorem \ref{pseudotorsor thm} will
apply to $X$ as well.
\end{remark}
\begin{lemma}\label{calculate moment}
The moment map for the action of $GL_n$ on $T^*\wt{\PS}_n(X,V)$
lifting the action on $\wt{\PS}_n(X,V)$ is the map
$$T^*\wt{\PS}_n(X,V)\to \on{Hom}(\Oo^n[-1],\Oo^n\ot \Omega)\simeq
\gl_n^*$$ assigning to $(Q,i,q,r)$ the composition
$$\Oo^n[-1]\xrightarrow{q[-1]} Q[-1] \xrightarrow{r} \on{Cone}(V\oplus
\Oo^n \xrightarrow{(i,q)}
Q)\ot \Omega \xrightarrow{\on{can}} \Oo^n\ot\Omega,$$ where $\on{can}$ is the canonical
map from the cone of a morphism to (a summand of) its source.
\end{lemma}

\begin{prop}\label{PB prop} The stack $\PB_n(X,V)$ is canonically identified with the
Hamiltonian reduction $$T^*\wt{\PS}_n(X,V)/\!\!/_{0}\, GL_n$$ of
$T^*\wt{\PS}_n(X,V)$ at the moment value $0$, i.e. to the cotangent
stack to $\PS_n(X,V)$.
\end{prop}

\begin{proof}
Maps $r$ as in \eqref{r} that
satisfy the vanishing moment condition define and are defined by
maps
$$\ol{a}: Q[-1]\to \on{Cone}(V \xrightarrow{i} Q)\ot
\Omega.$$
The proposition then follows from Theorem \ref{perverse thm}.
\end{proof}

\begin{lemma}
The canonical morphisms
\begin{equation}\label{exts}
\Ext^1(Q,Q\ot\Omega) \leftarrow \Ext^1(Q,\Oo^n\ot\Omega) \rightarrow
\Ext^1(\Oo^n,\Oo^n\ot\Omega) \end{equation} are both isomorphisms.
\end{lemma}

\begin{proof} First observe that all three vector spaces have dimension
$n^2$. The left arrow is surjective (hence an isomorphism) due to
the vanishing of $\Ext^2$'s. The right arrow is dual to the map
$\Hom(\Oo^n,Q) \leftarrow \Hom(\Oo^n,\Oo^n)$, which is injective
(hence an isomorphism) since the kernel of $\Oo^n\to Q$ has no
global sections.
\end{proof}

Let $\cJ^1$ denote the sheaf of one-jets of functions on $X$.  Let
$[X]\in H^1(\Omega^1)$ be the canonical nonzero class (i.e. the fundamental class of $X$).
\begin{lemma}\label{class is fundamental}
The image of the extension class
\bd
[\cJ^1(Q)] \overset{\on{def}}{=} [Q\otimes\cJ^1]\in
\on{Hom}(Q[-1],Q\otimes\Omega)
\ed
 in $\on{Ext}^1(\theo^n,\theo^n\otimes\Omega)$
under the isomorphism of \eqref{exts} equals $[X]\otimes\theo^n$.
\end{lemma}
\begin{proof}
In the case $Q=\theo_p$, $p\in X$ the assertion follows immediately
from the observation that the tautological extension of $\theo_p$ by
$\theo_p\ot \Omega$, Serre dual to the identity map of $\theo_p$, is
given by one-jets. This class is identified using \eqref{exts} (and
the functoriality of Serre duality) with the Serre dual to the
identity map of $\theo_X$, i.e. the fundamental class $[X]$. The
statement of the lemma is now immediate for $Q=\oplus_i \theo_{p_i}$
where $p_1,\dots, p_n$ are points of $X$, and follows for general
$Q$ by continuity.
\end{proof}

We will now perform hamiltonian reduction from $T^*\wt{\PS}_n(X,V)$
at a nontrivial moment value, corresponding to the dual of the
determinant character of $GL_n$. More precisely, the extension class
$$[X]^n:=[X]\otimes
\theo^n\in\on{Hom}(\Oo^n[-1],\Oo^n\ot\Omega)=\gl_n^*$$ corresponds
to the determinant character; we will take the hamiltonian reduction
of $T^*\wt{\PS}_n(X,V)$ at $-[X]^n.$

\begin{thm}\label{pseudotorsor thm} The Calogero-Moser stack $\CM_n(X,V)$ is isomorphic, as a
pseudotorsor over $\PB_n(X,V)=T^*\PS_n(X,V)$, to the hamiltonian
reduction
$$T^*\wt{\PS}_n(X,V)/\!\!/_{-[X]^n}\, GL_n,$$ i.e. to the twisted
cotangent bundle to $\PS_n(X,V)$ associated to the dual determinant
line bundle.
\end{thm}

\begin{proof}

Let $(Q,i,q)\in \wt{\PS}_n(X,V)$ be a framed torsion sheaf with a
basis. We wish to show the equivalence of CM data $a:Q\ot\D^1[-1]\to
\on{Cone}(V\to Q)$, satisfying the condition on its restriction to
zeroth order differential operators, and of cotangent data $r$ as
above, satisfying the moment condition.

Since $\cJ^1$ is the dual of $\D^1$, the data of the map $a$ is
equivalent to that of a map
$$a_1:Q[-1]\to \on{Cone}(\cJ^1 V \xrightarrow{\cJ^1(i)} \cJ^1 Q)$$
such that the composite map
\bd
Q[-1]\rightarrow \on{Cone}(\cJ^1\otimes (V\rightarrow Q))\rightarrow
\on{Cone}(V\rightarrow Q)
\ed
is the standard one (``inclusion of $Q$'').
The datum of $a_1$ is similarly equivalent to the choice of a map
$$a_2: Q[-1]\to \on{Cone}(V\ot \Omega
\xrightarrow{\cJ^1(i)|_{V\ot\Omega}} \cJ^1 Q)$$
such that the composite
\bd
a_2: Q[-1]\to \on{Cone}(V\ot \Omega
\xrightarrow{\cJ^1(i)|_{V\ot\Omega}} \cJ^1 Q)
\rightarrow Q[-1]
\ed
is the identity.

Next, note that our condition on $a_2$ implies that the composite
\begin{equation}\label{partial composite}
Q[-1]\xrightarrow{a_2}\on{Cone}\big(V\ot \Omega
\xrightarrow{\cJ^1(i)|_{V\ot\Omega}} \cJ^1 Q\big)\rightarrow
\on{Cone}\left(Q\otimes\Omega\rightarrow \cJ^1(Q)\right)
\end{equation}
is inverse to the canonical quasi-isomorphism
$\on{Cone}(Q\otimes\Omega\rightarrow \cJ^1(Q))\simeq Q[-1]$.
Consequently, if we compose \eqref{partial composite} with the further
projection to $Q\otimes\Omega$, the resulting class in
$\on{Hom}(Q[-1],Q\otimes\Omega)$ is exactly $[\cJ^1(Q)]$.  Let
 $\ol{a}_2:Q[-1]\to V\ot\Omega$ be the result of composing $a_2$ with the
projection on $V\otimes\Omega$.  The above then shows that the composite
\bd
Q[-1]\xrightarrow{\ol{a}_2}V\ot\Omega \xrightarrow{i\otimes\Omega}
 Q\otimes\Omega
\ed
also gives
$[\cJ^1(Q)]\in\on{Hom}(Q[-1],Q\ot\Omega)=\on{Ext}^1(Q,Q\ot\Omega)$.

We now define a map $$\wt{r}:Q[-1]\to (V\ot \Omega)\oplus (\Oo^n\ot
\Omega)$$ by $\wt{r} = (\ol{a}_2, -q^\vee)$, where $q^\vee\in
\Ext^1(Q,\Oo^n\ot\Omega)$ corresponds to $\cJ^1 Q$ under the first
identification of \eqref{exts}. It follows that the
composite
\bd
Q[-1]\xrightarrow{\wt{r}}(V\ot \Omega)\oplus (\Oo^n\ot
\Omega)\xrightarrow{(i+ q)\ot\Omega}
Q\otimes\Omega
\ed
 is
\bd (i\otimes\Omega)\circ\ol{a}_2 - (q\otimes\Omega)\circ q^\vee =
[\cJ^1(Q)] - [\cJ^1(Q)]=0. \ed As a result, we get a lift $r$ of
$\wt{r}$ to $\on{Cone}\big((i+ q)\ot\Omega\big)$.  To see that $r$
satisfies the moment condition, we compute that \bd \on{can}\circ
r\circ q[-1] = -q^\vee\circ q[-1], \ed which equals $-[X]^n$ by
Lemma \ref{class is fundamental}.

Conversely, given $r$ as above, we find that its restriction to a
map $Q[-1]\to V\ot\Omega$ composes with $q\ot \Omega$ to define the
one-jet map $Q[-1]\to Q\ot\Omega$. It follows that we may use it to
define a map $a_2$ as above, and hence a map $a$ satisfying the
restriction condition as desired.

Clearly all the identifications are $GL_n$-equivariant and preserve
the torsor structures over the cotangent bundle of $\PS_n(X,V)$. The
theorem follows.
\end{proof}

\section{Explicit Description of CM Matrices for ${\mathbf
A}^1$}\label{explicit}
 In this section we illustrate in detail how
our description of $\D$-bundles specializes, in the case of the
affine line, to the familiar Calogero-Moser space description. In
other words we want to spell out the correspondence between:
\begin{enumerate}
\item
quadruples
\bd
(X,Y,i,j)\in \gl_n\times\gl_n\times\C^n \times(\C^n)^*
\ed
satisfying
the CM relation $[X,Y] - ij + I = 0$ (which we call
``CM quadruples'') up to simultaneous conjugation
and
\item our data for $\Oo$-framed $\D$-bundles on ${\mathbf A}^1$ with
local second Chern class $n$,
 consisting of a length $n$ torsion sheaf $Q$, a framing map
$\theo\rightarrow Q$,
and an action map $\D^1\otimes Q[-1] \rightarrow \on{Cone}(\theo\rightarrow
Q)$ (defined in the derived category, i.e. up to homotopy) 
that agrees with the usual map on $\theo\otimes Q\subset
\D^1\otimes Q$.
\end{enumerate}
To summarize the relationship concisely:
 the matrix $X$ describes the sheaf $Q$,
the vector $i$ describes the framing map (justifying our general
notation for the latter), and $Y,j$ provide the components of the
action morphism $a$.
\begin{remark}
It is important to note that a CM quadruple corresponds to a unique triple $(Q,i,a)$
but where $a$ is only well defined {\em up to homotopy}.  What we will find below
is that, once we replace $Q[-1]$ by a resolution, we can produce a canonical representative
of the homotopy class of $(Q, i, a)$, and from this we can read off the CM quadruple 
completely explicitly.
\end{remark}

Choosing a basis for $H^0(Q)$, we may identify $Q$ with $\C^n$. A
choice of $\C[x]$-module structure on $\C^n$ is given by an $n\times
n$ matrix $X$, and we then get a presentation of $Q$ by
\begin{equation}\label{Q presentation}
0\rightarrow \C[x]^n\xrightarrow{x-X} \C[x]^n \rightarrow Q
\rightarrow 0
\end{equation}
where the last map takes $f(x)\otimes v$ to $f(X)v$ for $v\in\C^n
=Q$. Our Koszul data, then, take the form
\begin{equation}\label{Koszul diagram}
\xymatrix{\D^1\otimes\C[x]^n \ar[r]^{\ul{a}_0} \ar[d]^{1\otimes
(x-X)} &
\C[x] \ar[d]_{\ul{i}} \\
\D^1\otimes\C[x]^n\ar[r]^{\ul{a}_1} & \C^n.}
\end{equation}
Since all maps are required to be $\C[x]$-linear, the map $\ul{a}_0$
is determined by its values on vectors in $\C^n$ and elements of the
form $\partial\otimes v$ for $v\in\C^n$.  Since we require the
morphism of complexes to give the standard one on $\theo\otimes Q$,
we find that $\ul{a}_0(v) = 0$ for $v\in \C^n$; we write
$\ul{a}_0(\partial\otimes v) = j(x)v$ where $j(x):\C^n \rightarrow
\C[x]$ is a $\C$-linear map.  This completely determines $\ul{a}_0$.
Similarly, we may write $i = \ul{i}(1)$, which determines $\ul{i}$.
By \eqref{Q presentation}, $\ul{a}_1(f(x)\otimes v) = f(X)v$.  We
may write $\ul{a}_1(\partial\otimes v) = Yv$ for all $v\in\C^n$,
where $Y$ is an appropriate $n\times n$ matrix.  Applying
$\ul{a}_1\circ(1\otimes(x-X))$ and $\ul{i}\circ \ul{a}_0$ to an
element of the form $\partial\otimes v$, we may compute that the
square \eqref{Koszul diagram} commutes if and only if the following
equation is satisfied:
\begin{equation}
I+XY-YX = i\cdot j(X).
\end{equation}
Hence our Koszul data are determined by a choice of $(X,Y,i,j(x))$
satisfying this equation.  It is immediate that one has the
following converse: given a CM quadruple $(X,Y,i,j)$ satisfying the
CM relation, the diagram of complexes defined as above (with $j(x) =
j$) gives our $\D$-bundle data.

It remains to show
that, the diagram corresponding to a quadruple $(X,Y,i,j(x))$ is homotopic to
a unique diagram corresponding to a CM quadruple.  A homotopy of our complex
is a diagonal ``lower left to upper right'' map
$h:\D^1\otimes \C[x]^n\rightarrow \C[x]$.  We first choose this map to vanish
on elements $f(x)\otimes v$, $v\in\C^n$, and to take the form
$h(\partial\otimes v) = h_\partial(x)(v)$ for a linear map
$h_\partial(x):\C^n\rightarrow \C[x]$.  A computation then shows that, modifying
our original complex by the homotopy, the data $(X,Y,i,j(x))$ are replaced by
\begin{displaymath}
(X,Y+h_\partial(X), i, j(x) + xh_\partial(x) - h_\partial(x)\cdot X).
\end{displaymath}
In particular, one may check that, writing $j(x) = j + xj'$ where
$j': \C^n\rightarrow \C[x]$ is a linear map, one can solve the
equation $j(x) + xh_\partial-h_\partial(x)\cdot X = j$ for
$h_\partial(x)$. Consequently, up to homotopy, we can replace $(X,Y,
i,j(x))$ by a CM quadruple giving the same homotopy class of Koszul
data.

 To complete the proof, we start
from a CM quadruple $(X,Y,i,j)$ and consider an arbitrary homotopy
$h$.  As before, we write $h_\partial(x)(v)$ for $h(\partial\otimes
v)$. The formula above shows that, after modifying by the homotopy,
$j$ is replaced by $j+ xh_\partial(x)-h_\partial(x)\cdot X$, which,
in particular, has degree greater than zero in the variable $x$
unless $h_\partial=0$.  So if modifying $(X,Y,i,j)$ by $h$ gives
another CM quadruple $(X, Y', i, j')$ then $h_\partial=0$. Moreover,
it is easy to check that the new map $J'$ obtained after applying
$h$ satisfies $J'(v) = xh(v)-h(Xv)$. Writing the linear map $h:
\C^n\rightarrow \C[x]$ as a polynomial $\sum_{i=0}^k x^i h_i$ with
covector coefficients $h_i$, we find that if $h_k(v)\neq 0$ then
$\deg xh(v) =k+1 >\deg h(Xv)$. Since we require $J'(v)= 0$ for all
$v\in\C^n$, we conclude that $h(v) = 0$ for all $v\in\C^n$. Together
with the vanishing of $h_\partial$, this proves that $h=0$, i.e.
there is a unique CM quadruple in a given homotopy class of our
Koszul data. This completes the correspondence.

\end{document}